\documentclass[reqno]{amsart}
\usepackage{graphicx}

\begin{document}

\centerline{\LARGE\bf A Computer Approach to}

\medskip
\centerline{\LARGE\bf Determine the Densest Translative Tetrahedron Packings}

\bigskip\medskip
\centerline{\large Chuanming Zong}

\vspace{0.6cm}
\begin{minipage}{14cm}
{\bf Abstract.} In 1900, as a part of his 18th problem, Hilbert proposed the question to determine the densest congruent (or translative) packings of a given solid, such as the unit ball or the regular tetrahedron of unit edges. Up to now, our knowledge about this problem is still very limited, excepting the ball case. It is conjectured that, for some particular solids such as tetrahedra, cuboctahedra and octahedra, their maximal translative packing densities and their maximal lattice packing densities are identical. To attack this conjecture, this paper suggests a computer approach to determine the maximal local translative packing density of a given polytope, by studying associated color graphs and applying optimization. In particular, all the tetrahedral case, the cuboctahedral case and the octahedral case of the conjecture have been reduced into finite numbers of manageable optimization problems.
\end{minipage}

\vspace{1cm}
\centerline{\Large\bf 0. Introduction}

\vspace{0.4cm}\noindent
More than 23,00 years ago, Aristotle \cite{arist} claimed that identical regular tetrahedra can fill the whole space without gap. Unfortunately, this statement is wrong. In other words, no matter how to arrange the identical regular tetrahedra, they can not tile the whole space. This was discovered by M\"uller in fifteenth century (see Lagarias and Zong \cite{lazo}). In 1611, Kepler \cite{kepler} studied the sphere packing problem and conjectured that the density of the densest sphere packing is $\pi /\sqrt{18}$. In 1840, Gauss proved Kepler's conjecture for the lattice packings, by studying positive definite ternary quadratic forms. However, the final solution for Kepler's conjecture was discovered only at the end of the last century by Hales \cite{hales}.

In his 1900 ICM talk, David Hilbert \cite{hilbert} proposed 23 mathematical problems. At the end of his 18th problem, based on Aristotle's mistake on tetrahedron packings, Kepler's conjecture and Gauss' work on lattice sphere packings, he asked \lq\lq {\it How can one arrange most densely in space an infinite number of equal solids of given form, e.g., spheres with given radii or regular tetrahedra with given edges $($or in prescribed position$)$, that is, how can one so fit them together that the ratio of the filled to the unfilled space may be as great as possible}?"

\medskip
Let $K$ denote a {\it convex body} in the three-dimensional Euclidean space $\mathbb{E}^3$ containing the origin, with boundary $\partial (K)$, interior ${\rm int}(K)$ and volume ${\rm vol}(K)$, and let $P$ denote a centrally symmetric polytope centered at the origin. In particular, let $T$, $O$ and $C$ denote a regular tetrahedron, a regular octahedron and a regular cuboctahedron all with unit edges, respectively, and let $B$ denote the unit ball centered at the origin. Let $\delta^c(K)$, $\delta^t(K)$ and $\delta^l(K)$ denote the densities of the densest {\it congruent packings}, the densest {\it translative packings} and the densest {\it lattice packings} of $K$, respectively. It follows from their definitions that
$$\delta^l(K)\le \delta^t(K)\le \delta^c(K)\le 1\eqno (0.1)$$
holds for every convex body $K$. Moreover, both $\delta^l(K)$ and $\delta^t(K)$ are invariants under nonsingular affine linear transformations, while $\delta^c(K)$ for some $K$ is not. Then, Hilbert's problem can be restated as: {\it To determine the values of $\delta^c(K)$, $\delta^t(K)$ and $\delta^l(K)$ for a given solid $K,$ such as a sphere or a regular tetrahedron.}

The first approach to Hilbert's problem was made by Minkowski \cite{mink} in 1904. He defined
$$D(K)=\{ {\bf x}-{\bf y}:\ {\bf x}, {\bf y}\in K\}$$
to be the {\it difference body} of $K$ and proved
$$\delta^t(K)={{2^3{\rm vol}(K)}\over {{\rm vol}(D(K))}}\cdot \delta^t(D(K))\eqno (0.2)$$
and
$$\delta^l(K)={{2^3{\rm vol}(K)}\over {{\rm vol}(D(K))}}\cdot \delta^l(D(K)).\eqno (0.3)$$
Clearly, the difference body $D(K)$ is always centrally symmetric. Then, for a centrally symmetric convex body $K$ he discovered a criterion for its densest lattice packings. As an application, he determined the density of the densest lattice packings of an octahedron $O$. In other words, he proved
$$\delta^l(O)={{18}\over {19}}.\eqno (0.4)$$

On page 312 of \cite{mink}, Minkowski wrote \lq\lq {\it If $K$ is a tetrahedron, then ${1\over 2}D(K)$ is an octahedron with faces parallel to the faces of the tetrahedron.}"  By routine computations, one can get ${\rm vol}(T)=\sqrt{2}/{12}$ and ${\rm vol}(O)=\sqrt{2}/3$. Then, by (0.3) and (0.4) Minkowski \cite{mink} made a conclusion that $\delta^l(T)={9/{38}}.$ Unfortunately, Minkowski made a mistake, which was discovered by Groemer \cite{groe} in 1962. The difference body of a tetrahedron is not an octahedron, but a cuboctahedron. In fact, it was already known to Estermann and S\"uss in 1928 that
$${{{\rm vol}(D(T))}\over {{\rm vol}(T)}}=20.\eqno (0.5)$$

In 1970, Hoylman \cite{hoyl} applied Minkowski's criterion to a cuboctahedron $C$. By considering $38$ cases with respect to the possible positions of the three vectors of a basis, he proved that
$$\delta^l (C)={{45}\over {49}}, \eqno (0.6)$$
$$\delta^l(T)={{18}\over {49}},\eqno (0.7)$$
and the optimal lattice is unique up to certain equivalence. It is noteworthy that in the densest lattice cuboctahedron (tetrahedron) packing each cuboctahedron (tetrahedron) touches $14$ others.

Recent years, much progress has been made in the study of $\delta^c(T)$ by Conway, Torquato, Chen,  Engel, Glotzer, Kallus, Gravel, Elser, Jiao and etc. So far, the best known bounds are
$$0.856347\ldots \le\delta^c(T)\le 1-2.6\times 10^{-25}.$$

Perhaps, to determine the value of $\delta^t(K)$ is not as challenging as that for $\delta^c(K)$. However, since
$\delta^t(K)$ is invariant under nonsingular affine linear transformations, it is important. For $\delta^t(T)$ and
$\delta^t(C)$, by (0.1), (0.2), (0.5), (0.6) and (0.7) one can deduce
$${{45}\over {49}}\le \delta^t(C)\le 1$$
and
$${{18}\over {49}}\le \delta^t(T)\le {2\over 5}.$$

In 2014, by studying the shadow region in a packing, Zong \cite{zong} obtained the first nontrivial upper bounds for $\delta^t(C)$ and $\delta^t(T)$. Namely,
$$\delta^t(C)\le {{90\sqrt{10}}\over {95\sqrt{10}-4}}\approx 0.9601527\ldots $$
and
$$\delta^t(T)\le {{36\sqrt{10}}\over {95\sqrt{10}-4}}\approx 0.3840610\ldots .$$
Recently, these bounds were improved by Dostert, Guzman, de Oliveira Filho and Vallentin \cite{dgov} into $0.9364207\ldots $ and $0.3745683\ldots $, respectively, by a computational approach.

\smallskip
Based on the known results about $\delta^l(O)$, $\delta^t(C)$ and $\delta^t(T)$, it is reasonable to make the following conjecture:

\medskip\noindent
{\bf Conjecture A.}
$$\delta^t(O)={{18}\over {19}}, \qquad\delta^t(C)={{45}\over {49}} \qquad {\rm and}\qquad \delta^t(T)={{18}\over {49}}.$$

\medskip
To attack this conjecture, this paper suggests a computer approach to determine the maximal local translative packing density of a given polytope, by studying associated color graphs and applying optimization. In particular, all the tetrahedral case, the cuboctahedral case and the octahedral case of the conjecture have been reduced into finite numbers of manageable optimization problems.

\vspace{1cm}
\centerline{\Large\bf 1. Local Packing and Local Density}

\vspace{0.4cm}\noindent
By (0.2), without loss of generality, we may assume that $K$ is centrally symmetric. Furthermore, by {\it John's theorem}, we may assume that
$$B\subseteq K\subseteq \sqrt{3}B.\eqno(1.1)$$
Let $\mathcal{K}$ denote the space of all such three-dimensional centrally symmetric convex bodies associated with the {\it Hausdorff metric}.
It is well-known and easy to see that both $\delta^t(K)$ and $\delta^l(K)$ are continuous on $\mathcal{K}$.

\vspace{0.6cm}\noindent
{\large\bf 1.1. Local Density}

\medskip
\noindent
To study $\delta^t(K)$ for a given $K$, the most natural approach is localization. Let $X$ be a discrete set of points in $\mathbb{E}^3$ which contains the origin ${\bf o}$ such that $K+X$ is a packing, let $\Pi (X)$ denote the {\it Dirichlet-Voronoi cell} of ${\bf o}$ determined by $X$, and let $\mathfrak{X}$ denote the family of
all such sets $X$. In other words,
$$\Pi (X)=\{ {\bf y}:\ \| {\bf y}, {\bf o}\| \le \| {\bf y}, {\bf x}\|,\ {\bf x}\in X\}.$$
It is easy to see that $\Pi (X)$ is a convex polytope. Then, we define
$$\omega (K)=\min_{X\in \mathfrak{X}}{\rm vol}(\Pi (X))$$
and
$$\delta (K)={{{\rm vol}(K)}\over {\omega (K)}}.$$

It is easy to show that $\delta (K)$ is continuous on $\mathcal{K}$ as well, and
$$\delta^t(K)\le \delta (K)$$
holds for all centrally symmetric convex bodies $K$. However, $\delta (K)$ is no longer always invariant under nonsingular linear transformations.

In 1943, to study sphere packings, L. Fejes T\'oth \cite{fejes} made the following conjecture:

\medskip\noindent
{\bf The Dodecahedral Conjecture.} {\it In any packing of unit balls the volume of each Voronoi cell has volume at least that of a regular dodecahedron of inradius one.}

\medskip
In other words, he believed that
$$\omega (B)=5\sqrt{3}\left(3\tan^2(\pi /5)-1\right)$$
and therefore
$$\delta (B)={{4\pi }\over {15\sqrt{3}\left(3\tan^2(\pi /5)-1\right)}}.$$
This conjecture was proved by Hales and McLaughlin \cite{halesmc} only in 2010. Together with Hales' work on Kepler's conjecture, we have
$$\delta^t (B)\not= \delta(B).$$

Let $\mathcal{K}^*$ denote the subset of $\mathcal{K}$ consisting of all convex bodies $K$ satisfying $$\delta^t(K)<\delta (K).$$ Since both $\delta^t(K)$ and
$\delta (K)$ are continuous on $\mathcal{K}$, the subset $\mathcal{K}^*$ is open in $\mathcal{K}$ and therefore the complement set $\mathcal{K}\setminus \mathcal{K}^*$
is closed, if it is nonempty.

\medskip
\noindent
{\bf Remark 1.1.}  In the plane, it can be shown that
$$\delta^t(K)=\delta (K)$$
holds for circular discs (see Fejes T\'oth \cite{fejes}), regular hexagons and many others. On the other hand, it does not hold for narrow rectangles. For example, defining
$$Q_\epsilon =\left\{ (x, y):\ |x|\le \mbox{${1\over 2}$}\epsilon,\ |y|\le \mbox{${1\over {2\epsilon }}$}\right\},$$
we have ${\rm vol}(Q_\epsilon)=1$, $\delta^t(Q_\epsilon )=1$, and
$$\lim_{\epsilon\to 0}\delta (Q_\epsilon )=\infty .$$

In 1885, it was discovered by Fedorov \cite{fedo} that there are only four types of parallelohedra (we treat parallelotopes as particular hexagonal prisms, see Figure 1.1), which can translatively tile the whole three-dimensional space. For convenience, let $v(P)$, $e(P)$ and $f(P)$ denote the numbers of the vertices, the edges and the faces of a polytope $P$, respectively. Then we have the following facts about the parallelohedra:

\bigskip
\centerline{\begin{tabular}{|c|c|c|c|c|}\hline
& hexagonal & rhombic & elongated & truncated \\
$ P$ & prism & dodecahedron & octahedron & octahedron \\
\hline
$v(P)$ & 12 & 14 & 18 & 24 \\
\hline
$e(P)$ & 18 & 24 & 28 & 36 \\
\hline
$f(P)$ & 8 & 12 & 12 & 14 \\
\hline
\end{tabular}}

\bigskip\noindent
We note that, among the parallelohedra, all $v(P)$, $e(P)$ and $f(P)$ attain their maxima at truncated octahedra.

\begin{figure}[ht]
\centering
\includegraphics[height=9cm,width=9cm,angle=0]{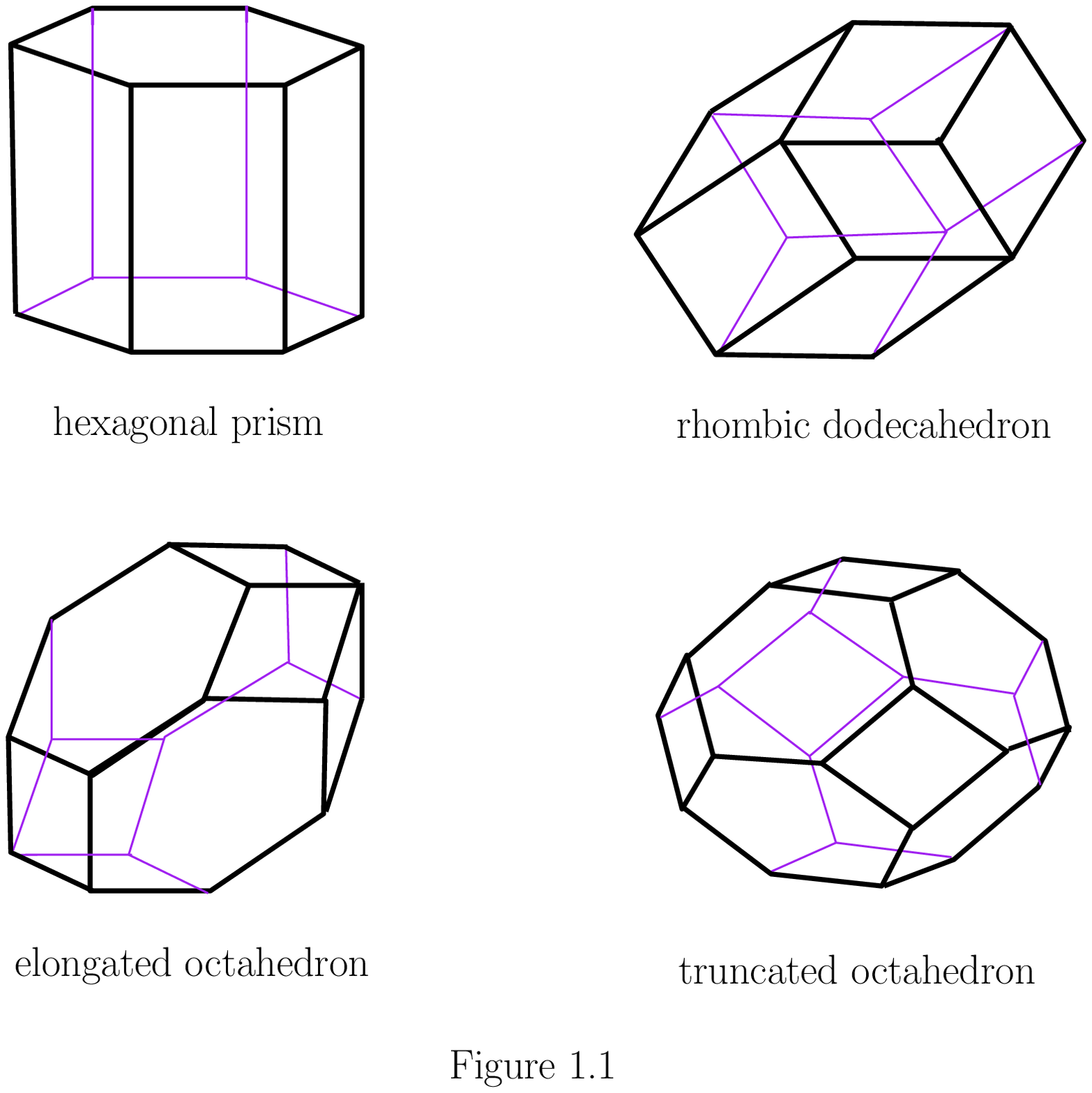}
\end{figure}

Lattice packings are particular translative packings. If $\Pi (\Lambda )$ is a regular truncated octahedron (with equal edges) whenever $K+\Lambda$ is one of the densest lattice packings of $K$,
by looking at their variations it is reasonable to believe that
$$\delta (K)=\delta^t(K)=\delta^l(K).$$

\noindent
{\bf Example 1.1 (Hoylman \cite{hoyl}).} Let $\Lambda_1 $ denote the lattice generated by ${\bf a}_1= \left( \mbox{${2\over 3}$}, 1, \mbox{${1\over 3}$}\right),$ ${\bf a}_2= \left( \mbox{${1\over 3}$}, \mbox{${2\over 3}$}, -1\right)$
and ${\bf a}_3= \left(\mbox{${4\over 3}$},  \mbox{${1\over 3}$}, -\mbox{${1\over 3}$}\right).$ It is well-known that $O+\Lambda_1$ is one of the densest lattice octahedron packings and $\Pi (\Lambda_1)$ is a truncated octahedron with equal edges. In fact, the Dirichlet-Voronoi cell $\Pi (\Lambda)$ will be a regular truncated octahedron whenever $O+\Lambda$ is one of the densest lattice packings of the regular octahedron $O$.

\medskip\noindent
{\bf Example 1.2 (Hoylman \cite{hoyl}).} Let $\Lambda_2$ denote the lattice generated by ${\bf a}_1= \left( 2, \mbox{${1\over 3}$}, \mbox{${1\over 3}$}\right),$  ${\bf a}_2= \left( \mbox{${2\over 3}$}, \mbox{${5\over 3}$}, \mbox{${5\over 3}$}\right)$ and ${\bf a}_3= \left( \mbox{${1\over 3}$}, -\mbox{${1\over 3}$}, 2\right).$ It is well-known that $C+\Lambda_2$ is one of the densest lattice cuboctahedron  packings and $\Pi (\Lambda_2)$ is also a truncated octahedron with equal edges. Similar to the octahedron case, the Dirichlet-Voronoi cell $\Pi (\Lambda )$ will be a regular truncated octahedron whenever $C+\Lambda$ is one of the densest lattice packings of the cuboctahedron $C$.

\medskip
These facts support the following conjecture:

\medskip\noindent
{\bf Conjecture Z.}
$$\delta (O)={{18}\over {19}} \qquad {\rm and}\qquad \delta (C)={{45}\over {49}}.$$

\bigskip
Clearly, Conjecture Z implies Conjecture A. Based on Graph Theory and Optimization, this paper develops a program to prove Conjecture Z.

\vspace{0.6cm}\noindent
{\large\bf 1.2. The Core of a Convex Body}

\medskip
\noindent
{\bf Definition 1.1.} Let $K$ be a centrally symmetric convex body centered at the origin. We define
$$K'=\{ {\bf x}:\ \langle {\bf x}, {\bf y}\rangle \le \langle {\bf y}, {\bf y}\rangle\ \mbox{for all ${\bf y}\in \partial (K)$}\}$$
and call it the core of $K$.

\medskip
In other words, let $H_{\bf x}$ denote the {\it bisector} of ${\bf o}$ and ${\bf x}$, and let $\overline{H_{\bf x}}$ denote the closed halfspace
bounded by $H_{\bf x}$ and containing the origin, then
$$K'=\bigcap_{{\bf y}\in \partial (K)}\overline{H_{2{\bf y}}}.$$

\smallskip
The core of a convex body has a number of simple properties. First of all, the core of a convex body is a convex body as well. Secondly,
$K'\subseteq K$ holds for every centrally symmetric convex body $K$. Thirdly, if $K_1\subseteq K_2$, for every ${\bf y}\in \partial (K_1)$
there is a $\gamma ({\bf y})\ge 1$ such that $\gamma ({\bf y}) {\bf y}\in \partial (K_2)$. In fact, this makes a {\it bijection} between $\partial (K_1)$ and $\partial (K_2)$.
Therefore, if
$$\langle {\bf x}, {\bf y}\rangle \le \langle {\bf y}, {\bf y}\rangle $$
holds for all ${\bf y}\in \partial (K_1)$, then
$$ \langle {\bf x}, \gamma ({\bf y}){\bf y}\rangle \le \langle \gamma ({\bf y}){\bf y}, \gamma ({\bf y}){\bf y}\rangle $$
holds for all $\gamma ({\bf y}){\bf y}\in \partial (K_2).$ Thus, it follows that
$$K_1'\subseteq K_2'.$$
Consequently, we have
$$B\subseteq K'\subseteq \sqrt{3}B, \eqno (1.2)$$
whenever $K$ satisfies (1.1).

Now, we present a characterization for $K'=K$.

\medskip
\noindent
{\bf Theorem 1.1.} {\it A centrally symmetric convex body $K$ is identical to its core $K'$ if and only if it is a ball centered at the origin.}

\medskip
\noindent
{\bf Proof.} It is apparent that $B'=B$. Now, we deal with the only if part. For convenience, for a point (vector) ${\bf w}$, let $B({\bf w})$ denote the ball with radius $\| {\bf w}\|$ and centered at ${\bf w}$.

Assume that $K$ is a convex body which is identical with its core $K'$. Let ${\bf x}$ be a fixed point on $\partial (K)$. For all ${\bf y}\in \partial (K)$, by the definition of the core we have
$$\langle {\bf x}, {\bf y}\rangle \le \langle {\bf y}, {\bf y}\rangle ,$$
$$\langle {\bf y}-{\bf x}, {\bf y}\rangle \ge 0, $$
$$\langle {\bf y}-{\bf x}, {\bf y}\rangle +\left\langle \mbox{$1\over 2$}{\bf x}, \mbox{$1\over 2$}{\bf x}\right\rangle \ge \left\langle \mbox{$1\over 2$}{\bf x}, \mbox{$1\over 2$}{\bf x}\right\rangle $$
and therefore
$$\left\langle {\bf y}-\mbox{$1\over 2$}{\bf x}, {\bf y}-\mbox{$1\over 2$}{\bf x}\right\rangle \ge \left\langle \mbox{$1\over 2$}{\bf x}, \mbox{$1\over 2$}{\bf x}\right\rangle .$$
In other words, the ball $B\left({1\over 2}{\bf x}\right)$ is contained in $K$. Therefore, the vector ${\bf x}$ is an norm of $K$ at the point ${\bf x}$.

Let $H$ be a two-dimensional plane passing the origin. Then, $\partial (K)\cap H$ is a centrally symmetric differentiable curve (there is a unique tangent line at every point) such that the vector ${\bf x}$ is an norm of $K\cap H$ at the point ${\bf x}$. Assume that ${\bf x}=(x,y)$ is defined by
$$\left\{
\begin{array}{ll}
x=r(\theta )\cos \theta ,&\\
\vspace{-0.35cm}
&\\
y=r(\theta )\sin \theta ,&
\end{array}\right.$$
where $r(\theta ) >0$ and $0\le \theta < 2\pi $. It can be deduced that $r(\theta )$ is differentiable. Then the tangent direction ${\bf v}=(x', y')$ of $K\cap H$ at ${\bf x}$ is defined by
$$\left\{
\begin{array}{ll}
x'=r'(\theta )\cos \theta -r(\theta )\sin \theta,&\\
\vspace{-0.35cm}
&\\
y'=r'(\theta )\sin \theta +r(\theta )\cos \theta.&
\end{array}\right.$$
Since the vector ${\bf x}$ is an norm of $K\cap H$ at the point ${\bf x}$, we have $\langle {\bf x}, {\bf v}\rangle =0$ and therefore
$$r(\theta )\cos \theta \bigl(r'(\theta )\cos \theta -r(\theta )\sin \theta\bigr) +r(\theta )\sin \theta \bigl(r'(\theta )\sin \theta +r(\theta )\cos \theta\bigr)=0,$$

$$r(\theta )r'(\theta )\left(\sin^2\theta +\cos^2\theta \right)=0$$
and finally
$$r'(\theta )=0.$$
This means that $K\cap H$ must be a circular domain for every $H$ and therefore $K$ must be a ball. The theorem is proved. \hfill{$\square$}

\medskip
Next, as examples, we determine the cores for the cube, the octahedron and the cuboctahedron, respectively. For convenience, we write
$$v=\max\{ |x|, |y|, |z|\}$$
and
$$w=\sqrt{x^2+y^2+z^2}.$$

\medskip
\noindent
{\bf Example 1.3.} Let $Q$ denote the cube $\{ (x,y,z):\ \max\{ |x|, |y|, |z|\}\le 1\}$. Then, by routine computation it can be deduced that
$$Q'=\left\{\mbox{${2\over {1+w}}$}(x,y,z):\ (x,y,z)\in Q\right\}.$$

\medskip
\noindent
{\bf Example 1.4.} Let $O$ denote the octahedron $\{ (x,y,z):\ |x|+|y|+|z|\le 1\}$. Then, one can deduce that
$$O'=\left\{ \mbox{${2\over {1+\sqrt{3}w}}$}(x,y,z):\ (x,y,z)\in O\right\}.$$

\medskip
\noindent
{\bf Example 1.5.} Let $C$ denote the cuboctahedron defined by
$$\{ (x,y,z):\ \max\{ |x|, |y|, |z|\}\le 1,\ |x|+|y|+|z|\le 2\}.$$ For ${\bf x}=(x,y,z)\in \partial (C)$,
let $\gamma ({\bf x})$ denote the positive number such that $\gamma ({\bf x}){\bf x}\in \partial (C')$.
Then, we have
$$ \gamma ({\bf x})=\left\{ \begin{array}{ll}
{{2}\over {w+1}},& \mbox{if ${\bf x}$ belongs to a square facet,}\\
\vspace{-0.35cm}
&\\
{4\over {2+\sqrt{3}w}}, & \mbox{if ${\bf x}$ belongs to a triangle facet and $(2-\sqrt{3})w+2v\le 2$,}\\
\vspace{-0.35cm}
&\\
{2\over {w+v}}, & \mbox{if ${\bf x}$ belongs to a triangle facet and $(2-\sqrt{3})w+2v\ge 2$.}
\end{array}\right.$$

\medskip
From the definitions of the core and the Dirichlet-Voronoi cell, one can easily deduce the following lemma, which will be useful to our program.

\medskip
\noindent
{\bf Lemma 1.1.} {\it If $K$ is a centrally symmetric convex body centered at the origin and $K+X$ is a packing, where ${\bf o}\in X$, then we have}
$$K'\subseteq \Pi (X).$$

\medskip\noindent
{\bf Proof.} Assume that $X=\{ {\bf o}, {\bf x}_1, {\bf x}_2, \ldots \}$. Since $K+X$ is a packing, for every point ${\bf x}_i\in X$ there is a corresponding positive number
$\alpha ({\bf x}_i)$ such that $\alpha ({\bf x}_i){\bf x}_i\in \partial (K)$ and $\alpha ({\bf x}_i)\le {1\over 2}$.

If ${\bf x}\in K'$, it follows by the definition of $K'$ that
$$\langle {\bf x}, \alpha ({\bf x}_i){\bf x}_i\rangle \le \langle \alpha ({\bf x}_i){\bf x}_i, \alpha ({\bf x}_i){\bf x}_i\rangle $$
holds for all ${\bf x}_i\in X$. Therefore, for all ${\bf x}_i\in X$, we have
$$\langle {\bf x}, {\bf x}_i\rangle \le \alpha ({\bf x}_i)\langle {\bf x}_i, {\bf x}_i\rangle \le \mbox{$1\over 2$}\langle {\bf x}_i, {\bf x}_i\rangle.$$
Consequently, we get
$$K'\subseteq \Pi (X).$$
The lemma is proved.\hfill{$\square$}

\vspace{0.6cm}\noindent
{\large\bf 1.3. Localization}

\medskip
\noindent
Let $r $ be a number with $r \ge 2$ and let $\mathfrak{X}_r$ denote the family of all the sets $X$ such that ${\bf o}\in X,$ $K+X$ is a packing and $X\subseteq r K$.
Then we define
$$\omega_r (K)=\min_{X\in \mathfrak{X}_r}{\rm vol}(\Pi (X)).$$

Clearly, for a given $K$ satisfying (1.1), $\omega_r (K)$ is a decreasing function of $r $. In fact, it easily follows from Lemma 1.1 that there is a
(smallest) number $\tau (K)\ge 2$ determined by $K$ such that
$$\omega_r (K)=\omega (K)$$
holds whenever $r \ge \tau (K)$. Furthermore, we define
$$\tau =\max_{K\in \mathcal{K}}\tau (K).$$

To study the local packing density $\delta (K)$, it is important to estimate $\tau (K)$. Of course, it would be particular interesting if one can determine the values
of $\tau $, $\tau (O)$ and $\tau (C)$.

\medskip
\noindent
{\bf Lemma 1.2 (Ball \cite{ball89}, Barthe \cite{barthe98}).} {\it If the $n$-dimensional unit ball is the ellipsoid of maximum volume inscribed in $K$, then
$${\rm vol} (K)\le 2^n,$$
where the equality holds if and only if $K$ is an $n$-dimensional cube of edge length two.}

\medskip
\noindent
{\bf Theorem 1.2.}
$$\tau \le 24.3.$$

\medskip\noindent
{\bf Proof.} For convenience, we write $\phi =24.3.$ Let $K$ be a fixed centrally symmetric convex body satisfying $B\subseteq K\subseteq \sqrt{3}B$ and consequently
$B\subseteq K'\subseteq \sqrt{3}B$. It was shown by Smith \cite{smit05} that
$$\delta (K)\ge \delta^l(K)\ge 0.53835\ldots .\eqno (1.3)$$

Assume that $X=\{ {\bf o}, {\bf x}_1, \ldots , {\bf x}_m\}$ is a discrete set of points in $\mathbb{E}^3$ such that $K+X$ is a packing, ${\bf x}_m\not\in \phi K$, and
$$\Pi (X)\neq \Pi (X\setminus \{{\bf x}_m\}).$$
Then, it follows that ${\bf x}_m\not\in \phi B$ and $\Pi (X)$ has a point on the bisector of ${\bf o}$ and ${\bf x}_m$. Consequently, the Dirichlet-Voronoi cell $\Pi (X)$
has a point ${\bf p}$ on the boundary of ${\phi\over 2}B$, and it contains the convex hull ${\rm conv}\{B,{\bf p}\}$ of $B$ and ${\bf p}$.

By routine computation, based on Lemma 1.2 and (1.3) we get
$${\rm vol}(K)\le 8,$$
$${\rm vol}(\Pi (X))>{\rm vol} ({\rm conv}\{ B, {\bf p}\})={\pi \over 3}\left( {\phi\over 2}+2+{2\over \phi }\right)$$
and
$${{{\rm vol}(K)}\over {{\rm vol}(\Pi (X))}}<{{24}\over {\pi \left({\phi\over 2}+2+{2\over \phi }\right)}}<0.53835 \le \delta (K).$$
This means that, for all considered $K$,
$$\omega_r (K)=\omega (K)$$
holds whenever $r \ge \phi.$ The theorem is shown.\hfill{$\square$}

\medskip
\noindent
{\bf Remark 1.2.} For the octahedron $O$ and the cuboctahedron $C$, by similar arguments it can be deduced that
$$\tau (O)\le 10\qquad {\rm and} \qquad \tau (C)\le 10.$$
Clearly, these upper bounds and the upper bound $24.3$ in Theorem 1.2 are too big. In the ball case, as it was proved by Hales and McLaughlin \cite{halesmc}, we have $\tau (B)=2$.
It would be rather surprising if there are examples in $\mathcal{K}$ for $\tau (K)>2$. Therefore it is reasonable to conjecture that
$$\tau=\tau (O)=\tau (C)=2.\eqno(1.4)$$

\medskip\noindent
{\bf Example 1.6.} The situation can be much different outside of $\mathcal{K}$. Let $\epsilon $ be a small positive number and define
$$P_\epsilon =\left\{ (x,y,z):\ |x|\le \mbox{${1\over {2\epsilon }}$},\ |y|\le \mbox{${1\over 2}$}, \ |z|\le \mbox{${1\over 2}$}\epsilon \right\}.$$
Clearly, we have ${\rm vol}(P_\epsilon )=1$ and $\delta^t (P_\epsilon )=\delta^l(P_\epsilon )=1$. On the other hand, by routine constructions, one can deduce that
$$\lim_{\epsilon \to 0}\omega (P_\epsilon )=0$$
and, for sufficiently small $\epsilon $,
$$\tau (P_\epsilon )>2.$$

\medskip
Next, we observe the localization from another view point. Let $\mathfrak{X}^m$ denote the family of sets $X$ such that ${\bf o}\in X$, $X$ has $m$ points and $K+X$ is a packing.
Similar to $\omega_\gamma (K)$, we may define
$$\varpi_m (K)=\min_{X\in \mathfrak{X}^m}{\rm vol}(\Pi (X)).$$

For a given $K$ satisfying (1.1), it is easy to see that $\varpi_m (K)$ is a decreasing sequence and there is a
(smallest) number $m(K)$ determined by $K$ such that
$$\varpi_m (K)=\omega (K)$$
holds whenever $m \ge m(K)$. Then we define
$$m^* =\max_{K\in \mathcal{K}}m(K).$$

Clearly, to determine the value of $m(K)$ is helpful to understand the local packings of $K$, and to determine the value of $m^*$ would be both interesting and important for understanding local packings in general.
By Theorem 1.2 and Remark 1.2, one can easily deduce the following result about $m(K)$, $m(O)$ and $m(C)$.

\medskip
\noindent
{\bf Corollary 1.1.} {\it For any fixed $K$ satisfying $(1.1)$, we have
$$m(K)\le 26^3.$$
In particular, for $O$ and $C$ we have
$$m(O)\le 11^3$$
and}
$$m(C)\le 11^3.$$

\medskip
\noindent
{\bf Remark 1.3.} Clearly, these upper bounds are far away from optimal. Similar to (1.4), it is reasonable to conjecture that
$$m^*=m(O)=m(C)=14.$$

\vspace{0.6cm}\noindent
{\large\bf 1.4. General Local Packings}

\medskip
\noindent
Let $P$ denote a three-dimensional centrally symmetric polytope with $2n$ faces $F_1$, $F_2$, $\ldots $, $F_{2n}$ and let ${\rm rint}(F_i)$ denote the relative interior of the face $F_i$. For convenience, we use ${\bf x}\prec F$ to abbreviate the fact that there is a positive number $\lambda $ such that $\lambda {\bf x}\in F$.

\medskip
\noindent
{\bf Definition 1.2.} {\it Let $X=\{ {\bf o}, {\bf x}_1, {\bf x}_2, \ldots , {\bf x}_m\}$ be a discrete set such that $P+X$ is a packing. We call $P+X$ a reduced local packing if $\Pi (X)$ is a polytope and
$$\Pi(X\setminus \{ {\bf x}_i\})\not= \Pi (X)$$
holds for every  ${\bf x}_i\in X\setminus \{{\bf o}\}$.}

\medskip
\noindent
{\bf Definition 1.3.} {\it A reduced local packing $P+X$ is called a general local packing if

\medskip\noindent
{\bf 1.} $\Pi (X)$ is a simple polytope $($three faces meet at each vertex$)$.

\noindent
{\bf 2.} For every ${\bf x}_i\in X\setminus \{{\bf o}\}$ we have
$${\bf x}_i\prec \bigcup_{k=1}^{2n}{\rm rint}(F_k).$$

\noindent
{\bf 3.} For every distinct pair $\{{\bf x}_i, {\bf x}_j\}$ in $X$ we have }
$${\bf x}_i-{\bf x}_j\prec \bigcup_{k=1}^{2n}{\rm rint}(F_k).$$

\medskip
\noindent
{\bf Theorem 1.3.} {\it For every reduced local packing $P+X$, where $X=\{ {\bf o}, {\bf x}_1, {\bf x}_2, \ldots , {\bf x}_m\}$, there is a sequence of general local packings
$P+Y_k$, where $Y_k=\{ {\bf o}, {\bf y}^k_1, {\bf y}^k_2, \ldots , {\bf y}^k_m\}$, such that
$$\lim_{k\to\infty }\| {\bf y}_i^k, {\bf x}_i\| =0$$
holds simultaneously for all $i=1, 2, \ldots, m.$}

\medskip
This theorem can be proved by a routine argument based on the fact that, for any positive small $\epsilon$ and $X_\epsilon =\{{\bf o}, (1+\epsilon ){\bf x}_1, (1+\epsilon ){\bf x}_2,
\ldots , (1+\epsilon ){\bf x}_m\}$, $P+X_\epsilon $ is also a reduced local packing and its translates are pairwise disjoint. Then, one can shift the $m$ translates
$P+(1+\epsilon ){\bf x}_i$ one by one in their small neighborhoods to produce a general local packing.

\medskip
This result guarantees the completeness of our method.

\vspace{1cm}
\centerline{\Large\bf 2. Color Graphs Associated with Local Packings}

\vspace{0.4cm}\noindent
First, let us recall two concepts which will be useful for this paper.

\medskip\noindent
{\bf Definition 2.1.} A planar graph $G$ is called a triangulated graph if adding any new edge will create crossing.

\medskip\noindent
{\bf Definition 2.2.} A triangulated graph $G$ is called a triangulated color graph if all its vertices and edges are colored by prescribed colors.

\vspace{0.6cm}\noindent
{\large\bf 2.1. Triangulated Color Graphs Associated to General Local Packings}

\medskip\noindent
Let $P$ be a three-dimensional centrally symmetric polytope with $2n$ faces $F_1$, $F_2$, $\ldots $, $F_{2n}$, where $F_{n+k}=-F_k$ for $k=1$, $2,$ $\ldots,$ $n$,
let $P+X$ be a general local packing, where $X=\{ {\bf o}, {\bf x}_1, {\bf x}_2,$ $\ldots,$ ${\bf x}_m\}$, and let $Q_i$ denote the face of $\Pi (X)$ which is on
the bisector of ${\bf o}$ and ${\bf x}_i$.

\medskip
Let $W=\{ {\bf w}_1, {\bf w}_2, \ldots , {\bf w}_n\}$ be a set of $n$ different colors. The general local packing $P+X$ determines a triangulated color graph $G$ as following: {\it

\begin{enumerate}
\item[\bf 1.] $G$ has $m$ distinct vertices ${\bf v}_1$, ${\bf v}_2$, $\ldots$, ${\bf v}_m$ (corresponding to ${\bf x}_1$, ${\bf x}_2$, $\ldots $, ${\bf x}_m$, respectively);
\item[\bf 2.] The vertex ${\bf v}_i$ is in color ${\bf w}_k$ if and only if ${\bf x}_i\prec F_k\cup F_{n+k};$
\item[\bf 3.] Two different vertices ${\bf v}_i$ and ${\bf v}_j$ is connected if and only if the two faces $Q_i$ and $Q_j$ of $\Pi (X)$ are connected;
\item[\bf 4.] If ${\bf v}_i$ and ${\bf v}_j$ are connected by an edge ${\bf e}_{i,j}$ (of course, ${\bf e}_{i,j}$ and ${\bf e}_{j,i}$ are the same),
we will color it by ${\bf w}_k$ if and only if ${\bf x}_i-{\bf x}_j\prec F_k\cup F_{n+k}.$
\end{enumerate}}

\noindent
We call $G$ a color graph associated to $P+X$.

\medskip\noindent
{\bf Example 2.1.} We take $$O=\{ (x, y, z):\ |x_1|+|y|+|z|\le 1 \}$$ and enumerate its faces by $F_1=\{ (x, y, z):\ x+y+z=1,\ x_i\ge 0\}$,
$F_2=\{ (x, y, z):\ -x+y+z=1,\ x\le 0,\ y \ge 0,\ z\ge 0\}$,
$F_3=\{ (x, y, z):\ x-y+z=1,\ x\ge 0,\ y \le 0,\ z\ge 0\}$,
$F_4=\{ (x, y, z):\ x+y-z=1,\ x\ge 0,\ y \ge 0,\ z\le 0\}$, and $F_{4+i}=-F_i$ for $i=1,$ $2,$ $3,$ $4$. Furthermore, we take $$X=\{ {\bf o}, {\bf x}_1, {\bf x}_2, \ldots, {\bf x}_{14}\},$$
where ${\bf x}_1=\left( {2\over 3}, 1, {1\over 3} \right),$
${\bf x}_2=\left( {4\over 3}, {1\over 3}, -{1\over 3}\right),$
${\bf x}_3=\left( {2\over 3}, -{2\over 3}, -{2\over 3}\right),$
${\bf x}_4=\left(1, -{1\over 3}, {2\over 3} \right),$
${\bf x}_5=\left( {1\over 3}, {1\over 3}, {4\over 3}\right), $
${\bf x}_6=\left( -{1\over 3}, -{2\over 3}, 1\right),$
${\bf x}_7=\left( {1\over 3}, -{4\over 3}, {1\over 3}\right),$
${\bf x}_8=\left( -{2\over 3}, -1, -{1\over 3}\right),$
${\bf x}_9=\left( -{1\over 3}, -{1\over 3}, -{4\over 3}\right),$
${\bf x}_{10}=\left( {1\over 3}, {2\over 3}, -1 \right),$
${\bf x}_{11}=\left( -{1\over 3}, {4\over 3}, -{1\over 3}\right),$
${\bf x}_{12}=\left( -1, {1\over 3},  -{2\over 3}\right),$
${\bf x}_{13}=\left( -{4\over 3}, -{1\over 3}, {1\over 3}\right),$
${\bf x}_{14}=\left( -{2\over 3}, {2\over 3}, {2\over 3}\right).$
It is easy to see that $O+X$ is a general local packing. Let $W=\{ {\bf w}_1, {\bf w}_2, {\bf w}_3, {\bf w}_4\}$ denote a set of colors, where ${\bf w}_1$ is red,
${\bf w}_2$ is green, ${\bf w}_3$ is blue, and ${\bf w}_4$ is black. Then the triangulated color graph $G$ determined by $O+X$ can be illustrated by Figure 2.1.

\medskip
\begin{figure}[ht]
\centering
\includegraphics[height=6cm,width=8cm,angle=0]{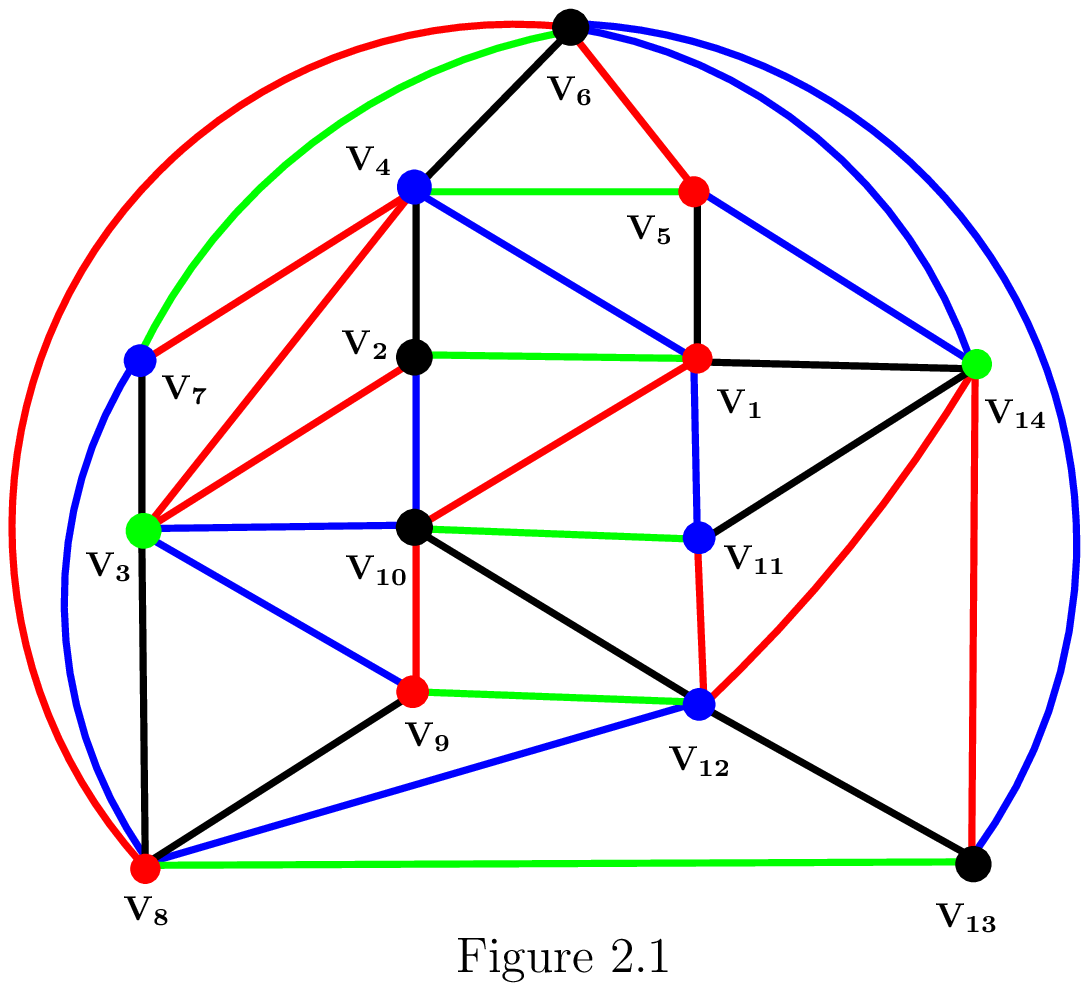}
\end{figure}

\medskip\noindent
{\bf Example 2.2.} We take $$C=\{ (x, y, z):\ |x_i|\le 1,\ |x|+|y|+|z|\le 2 \}$$ and enumerate its faces by
$F_1=\{ (x, y, z):\ x=1,\ |y|+|z|\le 1\}$,
$F_2=\{ (x, y, z):\ y=1,\ |x|+|z|\le 1\}$,
$F_3=\{ (x, y, z):\ z=1,\ |x|+|y|\le 1\}$,
$F_4=\{ (x, y, z):\ x_i\ge 0,\ x+y+z=2\}$,
$F_5=\{ (x, y, z):\ x\le 0,\ y\ge 0,\ z\ge 0,\ -x+y+z= 2\}$,
$F_6=\{ (x, y, z):\ x\ge 0,\ y\le 0,\ z\ge 0,\ x-y+z= 2\}$,
$F_7=\{ (x, y, z):\ x\ge 0,\ y\ge 0,\ z\le 0,\ x+y-z= 2\}$,
and $F_{7+i}=-F_i$ for $i=1,$ $2,$ $\ldots,$ $7$. Furthermore, we take $$X=\{ {\bf o}, {\bf x}_1, {\bf x}_2, \ldots, {\bf x}_{14}\},$$ where
${\bf x}_1=\left( 2, {1\over 3}, {1\over 3} \right),$
${\bf x}_2=\left( {2\over 3}, {5\over 3}, {5\over 3} \right),$
${\bf x}_3=\left( -{4\over 3}, {4\over 3}, {4\over 3} \right),$
${\bf x}_4=\left( {1\over 3}, -{1\over 3}, 2 \right),$
${\bf x}_5=\left( {5\over 3}, -{5\over 3}, {2\over 3} \right),$
${\bf x}_6=\left( -{1\over 3}, -2, {1\over 3}\right),$
${\bf x}_7=\left( -{5\over 3}, -{2\over 3}, {5\over 3} \right),$
${\bf x}_8=\left( -2, -{1\over 3}, -{1\over 3}\right),$
${\bf x}_9=\left( -{5\over 3}, {5\over 3}, -{2\over 3} \right),$
${\bf x}_{10}=\left( {1\over 3}, 2, -{1\over 3}\right),$
${\bf x}_{11}=\left( {5\over 3}, {2\over 3}, -{5\over 3} \right),$
${\bf x}_{12}=\left( -{1\over 3}, {1\over 3}, -2\right),$
${\bf x}_{13}=\left( -{2\over 3}, -{5\over 3}, -{5\over 3} \right),$
${\bf x}_{14}=\left( {4\over 3}, -{4\over 3}, - {4\over 3}\right).$
It is easy to see that $C+X$ is a general local packing. Let $W=\{ {\bf w}_1, {\bf w}_2, \ldots, {\bf w}_7\}$ denote a set of colors, where ${\bf w}_1$ is red,
${\bf w}_2$ is green, ${\bf w}_3$ is blue, ${\bf w}_4$ is yellow, ${\bf w}_5$ is purple, ${\bf w}_6$ is brown, and ${\bf w}_7$ is black.
Then the triangulated color graph $G$ determined by $C+X$ can be illustrated by Figure 2.2.

\medskip
\begin{figure}[ht]
\centering
\includegraphics[height=6cm,width=8cm,angle=0]{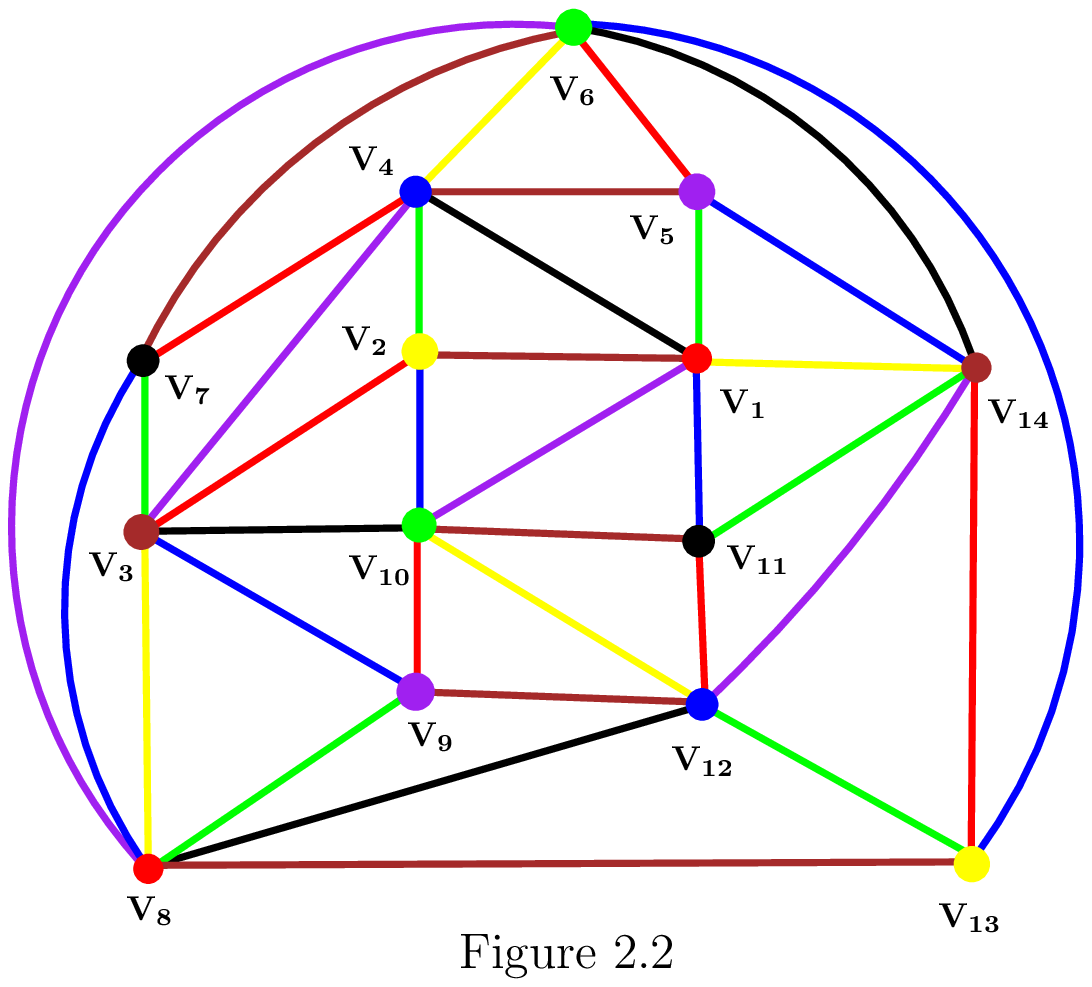}
\end{figure}

\vspace{0.6cm}\noindent
{\large\bf 2.2. Local Cells Associated to a Triangulated Graph}

\medskip
\noindent
The topological structure of a triangulated graph $G$ determines the geometric shape of its corresponding local cell $\Pi (X)$.

First, each vertex ${\bf v}_i$ of $G$ corresponds to a face $Q_i$ of the Voronoi cell $\Pi (X)$. In fact, the face $Q_i$ is on the bisector of ${\bf o}$ and ${\bf x}_i$,
where ${\bf x}_i$ is a point in $X$ corresponding to ${\bf v}_i$. Assume that ${\bf x}_i=(x_i, y_i, z_i)$ and write
$$d_i=\mbox{${1\over 2}$}\left(x_i^2+y_i^2+z_i^2\right),$$
then the points ${\bf x}=(x,y,z)$ of $Q_i$ satisfy the following equation
$$x_ix+y_iy+z_iz=d_i.\eqno(2.1)$$

Second, each edge ${\bf e}_{i,j}$ of $G$ corresponds to an edge $E_{i,j}$ of the Voronoi cell $\Pi (X)$. In fact, the points ${\bf x}=(x,y,z)$ of $E_{i,j}$ satisfy the following equations
$$\left\{
\begin{array}{l}
x_ix+y_iy+z_iz =d_i,\\
x_jx+y_jy+z_jz =d_j.
\end{array}
\right.$$

Third, each triangular face $\Delta_{i,j,k}$ of $G$ corresponds to a vertex ${\bf p}_{i,j,k}$ of the Voronoi cell $\Pi (X)$. In fact, the coordinates of the point ${\bf p}_{i,j,k}=(x_{i,j,k},y_{i,j,k},z_{i,j,k})$ are the solutions of
$$\left\{
\begin{array}{l}
x_ix+y_iy+z_iz =d_i,\\
x_jx+y_jy+z_jz =d_j,\\
x_kx+y_ky+z_kz =d_k.
\end{array}
\right.$$
That is
$$x_{i,j,k}={{a_{i,j,k}}\over {d_{i,j,k}}},\hspace{0.5cm}
y_{i,j,k}={{b_{i,j,k}}\over {d_{i,j,k}}},\hspace{0.5cm}
z_{i,j,k}={{c_{i,j,k}}\over {d_{i,j,k}}},\eqno(2.2)$$
where
$$ a_{i,j,k}=\left|
\begin{array}{ccc}
d_i & y_i & z_i\\
d_j & y_j & z_j\\
d_k & y_k & z_k
\end{array}
\right| ,\hspace{1cm}b_{i,j,k}=\left|
\begin{array}{ccc}
x_i & d_i & z_i\\
x_j & d_j & z_j\\
x_k & d_k & z_k
\end{array}
\right| ,$$

$$ c_{i,j,k}=\left|
\begin{array}{ccc}
x_i & y_i & d_i\\
x_j & y_j & d_j\\
x_k & y_k & d_k
\end{array}
\right| , \hspace{1cm} d_{i,j,k}=\left|
\begin{array}{ccc}
x_i & y_i & z_i\\
x_j & y_j & z_j\\
x_k & y_k & z_k
\end{array}
\right| .$$

Assume that there are $n(i)$ triangles taking ${\bf v}_i$ as a vertex in $G$. For convenience, we enumerate them as $\Delta_i^1$, $\Delta_i^2$, $\ldots$, $\Delta_i^{n(i)}$ in an anti-clock order,
and rewrite their corresponding vertices of $\Pi (X)$ as ${\bf p}_i^1=\bigl(x_i^1, y_i^1, z_i^1\bigr)$, ${\bf p}_i^2=\bigl(x_i^2, y_i^2, z_i^2\bigr)$, $\ldots $, ${\bf p}_i^{n(i)}=\bigl(x_i^{n(i)}, y_i^{n(i)}, z_i^{n(i)}\bigr)$,
respectively. It is easy to see that $Q_i$ is a convex $n(i)$-gon with vertices ${\bf p}_i^1$, ${\bf p}_i^2$, $\ldots $, ${\bf p}_i^{n(i)}$ in an anti-clock circular order.

Let us triangulate $Q_i$ into $n(i)-2$ triangles $\nabla_i^1={\rm conv}\bigl\{ {\bf p}_i^1, {\bf p}_i^2, {\bf p}_i^3\bigr\},$ $\nabla_i^2={\rm conv}\bigl\{ {\bf p}_i^1, {\bf p}_i^3, {\bf p}_i^4\bigr\},$ $\ldots $,
$\nabla_i^{n(i)-2}={\rm conv}\bigl\{ {\bf p}_i^1, {\bf p}_i^{n(i)-1}, {\bf p}_i^{n(i)}\bigr\}.$ Let $\| {\bf x}, {\bf y}\|$ denote the distance between two points ${\bf x}$ and ${\bf y}$, let $\langle {\bf x}, {\bf y}\rangle$
denote the inner product of the two vectors ${\bf x}$ and ${\bf y}$, and let $\mu (Q)$ denote the area of a polygon $Q$.
By routine computations, one can deduce that
$$\mu \bigl(\nabla_i^j\bigr)= {1\over 2} \sqrt{ \bigl\| {\bf p}_i^1, {\bf p}_i^{j+1}\bigr\|^2\cdot \bigl\| {\bf p}_i^1, {\bf p}_i^{j+2}\bigr\|^2-\bigl\langle {\bf p}_i^{j+2}-{\bf p}_i^1,
{\bf p}_i^{j+1}-{\bf p}_i^1\bigr\rangle^2},$$
where
$$\bigl\| {\bf p}_i^1, {\bf p}_i^j\bigr\|^2=\bigl( x_i^j-x_i^1\bigr)^2+\bigl( y_i^j-y_i^1\bigr)^2+\bigl( z_i^j-z_i^1\bigr)^2$$
and
\begin{align*}
\bigl\langle {\bf p}_i^{j+2}-{\bf p}_i^1, {\bf p}_i^{j+1}-{\bf p}_i^1\bigr\rangle = & \bigl( x_i^{j+2}-x_i^1\bigr)\bigl( x_i^{j+1}-x_i^1\bigr)+\bigl( y_i^{j+2}-y_i^1\bigr)\bigl( y_i^{j+1}-y_i^1\bigr)\\
& + \bigl( z_i^{j+2}-z_i^1\bigr)\bigl( z_i^{j+1}-z_i^1\bigr).
\end{align*}
Therefore, we have
$$\mu \left(Q_i\right)=\sum_{j=1}^{n(i)} \mu \bigl(\nabla_i^j\bigr)= {1\over 2} \sum_{j=1}^{n(i)}\sqrt{ \bigl\| {\bf p}_i^1, {\bf p}_i^{j+1}\bigr\|^2\cdot \bigl\| {\bf p}_i^1, {\bf p}_i^{j+2}\bigr\|^2-\bigl\langle {\bf p}_i^{j+2}-{\bf p}_i^1,{\bf p}_i^{j+1}-{\bf p}_i^1\bigr\rangle^2}.\eqno (2.3)$$

Divide $\Pi (X)$ into $m$ cones $P_i$ with a common vertex ${\bf o}$ and $m$ polygonal bases $Q_i$, $i=1$, $2,$ $\ldots $ , $m$, respectively. Clearly, the height $h_i$ of $P_i$ over its base $Q_i$ is ${1\over 2}\| {\bf o},
{\bf x}_i\|.$ In other words, for ${\bf x}_i=(x_i, y_i, z_i)$ we have
$$h_i={1\over 2}\sqrt{ x_i^2+y_i^2+z_i^2}.$$
Then one can deduce that
$${\rm vol}(P_i)={1\over {12}}\sum_{j=1}^{n(i)} \sqrt{\bigl(x_i^2+y_i^2+z_i^2\bigr)\Bigl( \bigl\| {\bf p}_i^1, {\bf p}_i^{j+1}\bigr\|^2\cdot \bigl\| {\bf p}_i^1, {\bf p}_i^{j+2}\bigr\|^2-\bigl\langle {\bf p}_i^{j+2}-{\bf p}_i^1,{\bf p}_i^{j+1}-{\bf p}_i^1\bigr\rangle^2 \Bigr) }$$
and therefore
$${\rm vol}(\Pi (X))={1\over {12}}\sum_{i=1}^m\sum_{j=1}^{n(i)} \sqrt{\bigl(x_i^2+y_i^2+z_i^2\bigr)\Bigl( \bigl\| {\bf p}_i^1, {\bf p}_i^{j+1}\bigr\|^2\cdot \bigl\| {\bf p}_i^1, {\bf p}_i^{j+2}\bigr\|^2-\bigl\langle {\bf p}_i^{j+2}-{\bf p}_i^1,{\bf p}_i^{j+1}-{\bf p}_i^1\bigr\rangle^2 \Bigr) }.$$

\vspace{0.6cm}\noindent
{\large\bf 2.3. Volume Optimization for Local Cells}

\medskip
\noindent
The colors of the triangulated color graph $G$ contain a lot of geometric information about its corresponding local cell $\Pi (X)$.

We recall that $P$ is a three-dimensional centrally symmetric convex polytope with $n$ pairs of faces $\pm F_1$, $\pm F_2$, $\ldots $, $\pm F_n$,
$P+X$ is a general local packing with $X=\left\{ {\bf o}, {\bf x}_1, {\bf x}_2, \ldots , {\bf x}_m\right\},$ $\Pi (X)$ is the local cell of ${\bf o}$ with respect to
$X$, $W=\left\{ {\bf w}_1, {\bf w}_2, \ldots , {\bf w}_n\right\}$ is a set of $n$ different colors associated to the $n$ pairs of faces of $P$,
and $G$ is the triangulated color graph associated to $P+X$.

Clearly, for the fixed polytope $P$, there are infinitely many general local packings $P+X$ having the same triangulated color graph $G$ and therefore having the same
volume formula for their local cells.

Assume that the equation of $F_i$ is
$$a_ix+b_iy+c_iz=\tau_i,$$
where $a_i$, $b_i$, $c_i$ and $\tau_i>0$ are suitable constants. Then, a vertex ${\bf v}_j$ of $G$ is in color ${\bf w}_i$ implies that ${\bf x}_j=(x_j, y_j,z_j)$
satisfies
$$\left| a_ix_j+b_iy_j+c_iz_j\right|\ge 2\tau_i\eqno(2.4)$$
and an edge ${\bf e}_{j,k}$ of $G$ is in color ${\bf w}_i$ implies that ${\bf x}_j=(x_j, y_j,z_j)$ and ${\bf x}_k=(x_k, y_k,z_k)$
satisfy
$$\left| a_i(x_j-x_k)+b_i(y_j-y_k)+c_i(z_j-z_k)\right|\ge 2\tau_i.\eqno(2.5)$$

Then, if $G$ has $\ell$ edges, to minimize the volume of $\Pi (X)$ among all general local packings $P+X$ associated with the same triangulated color graph $G$ is to minimize
$$f(X)={1\over {12}}\sum_{i=1}^m\sum_{j=1}^{n(i)} \sqrt{\bigl(x_i^2+y_i^2+z_i^2\bigr)\Bigl( \bigl\| {\bf p}_i^1, {\bf p}_i^{j+1}\bigr\|^2\cdot \bigl\| {\bf p}_i^1, {\bf p}_i^{j+2}\bigr\|^2-\bigl\langle {\bf p}_i^{j+2}-{\bf p}_i^1,{\bf p}_i^{j+1}-{\bf p}_i^1\bigr\rangle^2 \Bigr) }\eqno(2.6)$$
under $m$ constraints of (2.4) type and $\ell$ constraints of (2.5) type.

\medskip
\noindent
{\bf Remark 2.1.} For computer, possible alternative expressions for (2.4) and (2.5) are
$$\left( a_ix_j+b_iy_j+c_iz_j\right)^2\ge 4\tau_i^2$$
and
$$\bigl( a_i\left(x_j-x_k\right)+b_i\left(y_j-y_k\right)+c_i\left(z_j-z_k\right)\bigr)^2\ge 4\tau_i^2.$$

\medskip\noindent
{\bf Experiment 2.1.} Assume that $O+X$ is a general local octahedral packing associated with the color graph $G$ in Figure 2.1. Then $\Pi (X)$ has $14$ faces $Q_i$ determined by (2.1)
and $24$ vertices ${\bf p}_{i,j,k}=(x_{i,j,k}, y_{i,j,k}, z_{i,j,k})$ determined by (2.2), where $\{ i, j, k\}$ runs over the $24$ triples $\{1, 2, 4\}$, $\{1, 4, 5\}$, $\{1, 5, 14\}$,
$\{1, 2, 10\}$, $\{1, 10, 11\}$, $\{1, 11, 14\}$, $\{2, 3, 4\}$, $\{2, 3, 10\}$, $\{3, 4, 7\}$, $\{3, 7, 8\}$, $\{3, 8, 9\}$, $\{3, 9, 10\}$, $\{4, 5, 6\}$, $\{4, 6, 7\}$, $\{5, 6, 14\}$,
$\{6, 7, 8\}$, $\{6, 8, 13\}$, $\{6, 13, 14\}$, $\{8, 9, 12\}$, $\{8, 12, 13\}$, $\{9, 10,$ $12\}$, $\{10, 11, 12\}$, $\{11, 12, 14\}$ and $\{12, 13, 14\}$. Then, we rewrite the vertices of $Q_i$ as
$$\begin{array}{lllll}
Q_1:& {\bf p}_1^1={\bf p}_{1,2,10},&  {\bf p}_1^2={\bf p}_{1,10,11},&  {\bf p}_1^3={\bf p}_{1,11,14},&  {\bf p}_1^4={\bf p}_{1,5,14},\\
\vspace{-0.36cm}
&&&&\\
& {\bf p}_1^5={\bf p}_{1,4,5},& {\bf p}_1^6={\bf p}_{1,2,4};& & \\
\vspace{-0.36cm}
&&&&\\
Q_2:& {\bf p}_2^1={\bf p}_{1,2,4},&  {\bf p}_2^2={\bf p}_{2,3,4},&  {\bf p}_2^3={\bf p}_{2,3,10},&  {\bf p}_2^4={\bf p}_{1,2,10};\\
\vspace{-0.36cm}
&&&&\\
Q_3:& {\bf p}_3^1={\bf p}_{2,3,4},&  {\bf p}_3^2={\bf p}_{3,4,7},&  {\bf p}_3^3={\bf p}_{3,7,8},&  {\bf p}_3^4={\bf p}_{3,8,9},\\
\vspace{-0.36cm}
&&&&\\
& {\bf p}_3^5={\bf p}_{3,9,10},& {\bf p}_3^6={\bf p}_{2,3,10}; & &\\
\vspace{-0.36cm}
&&&&\\
Q_4:& {\bf p}_4^1={\bf p}_{1,2,4},&  {\bf p}_4^2={\bf p}_{1,4,5},&  {\bf p}_4^3={\bf p}_{4,5,6},&  {\bf p}_4^4={\bf p}_{4,6,7},\\
\vspace{-0.36cm}
&&&&\\
& {\bf p}_4^5={\bf p}_{3,4,7},& {\bf p}_4^6={\bf p}_{2,3,4}; & &\\
\vspace{-0.36cm}
&&&&\\
Q_5:& {\bf p}_5^1={\bf p}_{1,4,5},&  {\bf p}_5^2={\bf p}_{1,5,14},&  {\bf p}_5^3={\bf p}_{5,6,14},&  {\bf p}_5^4={\bf p}_{4,5,6};\\
\vspace{-0.36cm}
&&&&\\
Q_6:& {\bf p}_6^1={\bf p}_{4,5,6},&  {\bf p}_6^2={\bf p}_{5,6,14},&  {\bf p}_6^3={\bf p}_{6,13,14},&  {\bf p}_6^4={\bf p}_{6,8,13},\\
\vspace{-0.36cm}
&&&&\\
& {\bf p}_6^5={\bf p}_{6,7,8},& {\bf p}_6^6={\bf p}_{4,6,7}; & &\\
\vspace{-0.36cm}
&&&&\\
Q_7:& {\bf p}_7^1={\bf p}_{3,4,7},&  {\bf p}_7^2={\bf p}_{4,6,7},&  {\bf p}_7^3={\bf p}_{6,7,8},&  {\bf p}_7^4={\bf p}_{3,7,8};\\
\vspace{-0.36cm}
&&&&\\
Q_8:& {\bf p}_8^1={\bf p}_{3,7,8},&  {\bf p}_8^2={\bf p}_{6,7,8},&  {\bf p}_8^3={\bf p}_{6,8,13},&  {\bf p}_8^4={\bf p}_{8,12,13},\\
\vspace{-0.36cm}
&&&&\\
& {\bf p}_8^5={\bf p}_{8,9,12},& {\bf p}_8^6={\bf p}_{3,8,9}; & &\\
\vspace{-0.36cm}
&&&&\\
Q_9:& {\bf p}_9^1={\bf p}_{3,8,9},&  {\bf p}_9^2={\bf p}_{8,9,12},&  {\bf p}_9^3={\bf p}_{9,10,12},&  {\bf p}_9^4={\bf p}_{3,9,10};\\
\vspace{-0.36cm}
&&&&\\
Q_{10}:& {\bf p}_{10}^1={\bf p}_{1,2,10},&  {\bf p}_{10}^2={\bf p}_{2,3,10},&  {\bf p}_{10}^3={\bf p}_{3,9,10},&  {\bf p}_{10}^4={\bf p}_{9,10,12},\\
\vspace{-0.36cm}
&&&&\\
& {\bf p}_{10}^5={\bf p}_{10,11,12},& {\bf p}_{10}^6={\bf p}_{1,10,11}; & &\\
\vspace{-0.36cm}
&&&&\\
Q_{11}:& {\bf p}_{11}^1={\bf p}_{1,10,11},&  {\bf p}_{11}^2={\bf p}_{10,11,12},&  {\bf p}_{11}^3={\bf p}_{11,12,14},&  {\bf p}_{11}^4={\bf p}_{1,11,14};\\
\vspace{-0.36cm}
&&&&\\
Q_{12}:& {\bf p}_{12}^1={\bf p}_{8,9,12},&  {\bf p}_{12}^2={\bf p}_{8,12,13},&  {\bf p}_{12}^3={\bf p}_{12,13,14},&  {\bf p}_{12}^4={\bf p}_{11,12,14},\\
\vspace{-0.36cm}
&&&&\\
& {\bf p}_{12}^5={\bf p}_{10,11,12},& {\bf p}_{12}^6={\bf p}_{9,10,12}; & &\\
\vspace{-0.36cm}
&&&&\\
Q_{13}:& {\bf p}_{13}^1={\bf p}_{6,8,13},&  {\bf p}_{13}^2={\bf p}_{6,13,14},&  {\bf p}_{13}^3={\bf p}_{12,13,14},&  {\bf p}_{13}^4={\bf p}_{8,12,13};\\
\vspace{-0.36cm}
&&&&\\
Q_{14}:& {\bf p}_{14}^1={\bf p}_{1,5,14},&  {\bf p}_{14}^2={\bf p}_{1,11,14},&  {\bf p}_{14}^3={\bf p}_{11,12,14},&  {\bf p}_{14}^4={\bf p}_{12,13,14},\\
\vspace{-0.36cm}
&&&&\\
& {\bf p}_{14}^5={\bf p}_{6,13,14},& {\bf p}_{14}^6={\bf p}_{5,6,14}. & &
\end{array}$$

In this case, we have
$$n(1)=n(3)=n(4)=n(6)=n(8)=n(10)=n(12)=n(14)=6$$
and
$$n(2)=n(5)=n(7)=n(9)=n(11)=n(13)=4.$$
Then the volume of $\Pi (X)$ is
$$f_1(X)={1\over {12}}\sum_{i=1}^{14}\sum_{j=1}^{n(i)} \sqrt{\bigl(x_i^2+y_i^2+z_i^2\bigr)\Bigl( \bigl\| {\bf p}_i^1, {\bf p}_i^{j+1}\bigr\|^2\cdot \bigl\| {\bf p}_i^1, {\bf p}_i^{j+2}\bigr\|^2-\bigl\langle {\bf p}_i^{j+2}-{\bf p}_i^1,{\bf p}_i^{j+1}-{\bf p}_i^1\bigr\rangle^2 \Bigr) }.\eqno (2.7)$$
In addition, the colors of the vertices of $G$ imply
$$\left\{\begin{array}{l}
\left| x_1+y_1+z_1\right| \ge 2, \\
\left| x_2+y_2-z_2\right| \ge 2, \\
\left| x_3-y_3-z_3\right| \ge 2, \\
\left| x_4-y_4+z_4\right| \ge 2, \\
\left| x_5+y_5+z_5\right| \ge 2, \\
\left| x_6+y_6-z_6\right| \ge 2, \\
\left| x_7-y_7+z_7\right| \ge 2, \\
\left| x_8+y_8+z_8\right| \ge 2, \\
\left| x_9+y_9+z_9\right| \ge 2, \\
\left| x_{10}+y_{10}-z_{10}\right| \ge 2,\\
\left| x_{11}-y_{11}+z_{11}\right| \ge 2,\\
\left| x_{12}-y_{12}+z_{12}\right| \ge 2,\\
\left| x_{13}+y_{13}-z_{13}\right| \ge 2,\\
\left| x_{14}-y_{14}-z_{14}\right| \ge 2,
\end{array}
\right.\eqno(2.8)$$
and the colors of the edges of $G$ imply
$$\left\{\begin{array}{l}
\left| x_1-x_2-y_1+y_2-z_1+z_2\right| \ge 2, \\
\left| x_1-x_4-y_1+y_4+z_1-z_4\right| \ge 2,  \\
\left| x_1-x_5+y_1-y_5-z_1+z_5\right| \ge 2, \\
\left| x_1-x_{10}+y_1-y_{10}+z_1-z_{10}\right| \ge 2, \\
\left| x_1-x_{11}-y_1+y_{11}+z_1-z_{11}\right| \ge 2,  \\
\left| x_1-x_{14}+y_1-y_{14}-z_1+z_{14}\right| \ge 2,  \\
\left| x_2-x_3+y_2-y_3+z_2-z_3\right| \ge 2, \\
\left| x_2-x_4+y_2-y_4-z_2+z_4\right| \ge 2, \\
\left| x_2-x_{10}-y_2+y_{10}+z_2-z_{10}\right| \ge 2, \\
\left| x_3-x_4+y_3-y_4+z_3-z_4\right| \ge 2, \\
\left| x_3-x_7+y_3-y_7-z_3+z_7\right| \ge 2, \\
\left| x_3-x_8+y_3-y_8-z_3+z_8\right| \ge 2, \\
\left| x_3-x_9-y_3+y_9+z_3-z_9\right| \ge 2, \\
\left| x_3-x_{10}-y_3+y_{10}+z_3-z_{10}\right| \ge 2, \\
\left| x_4-x_5-y_4+y_5-z_4+z_5\right| \ge 2, \\
\left| x_4-x_6+y_4-y_6-z_4+z_6\right| \ge 2, \\
\left| x_4-x_7+y_4-y_7+z_4-z_7\right| \ge 2, \\
\left| x_5-x_6+y_5-y_6+z_5-z_6\right| \ge 2, \\
\left| x_5-x_{14}-y_5+y_{14}+z_5-z_{14}\right| \ge 2, \\
\left| x_6-x_7-y_6+y_7-z_6+z_7\right| \ge 2, \\
\left| x_6-x_8+y_6-y_8+z_6-z_8\right| \ge 2, \\
\left| x_6-x_{13}-y_6+y_{13}+z_6-z_{13}\right| \ge 2, \\
\left| x_6-x_{14}-y_6+y_{14}+z_6-z_{14}\right| \ge 2, \\
\left| x_7-x_8-y_7+y_8+z_7-z_8\right| \ge 2, \\
\left| x_8-x_9+y_8-y_9-z_8+z_9\right| \ge 2, \\
\left| x_8-x_{12}-y_8+y_{12}+z_8-z_{12}\right| \ge 2, \\
\left| x_8-x_{13}-y_8+y_{13}-z_8+z_{13}\right| \ge 2, \\
\left| x_9-x_{10}+y_9-y_{10}+z_9-z_{10}\right| \ge 2, \\
\left| x_9-x_{12}-y_9+y_{12}-z_9+z_{12}\right| \ge 2, \\
\left| x_{10}-x_{11}-y_{10}+y_{11}-z_{10}+z_{11}\right| \ge 2, \\
\left| x_{10}-x_{12}+y_{10}-y_{12}-z_{10}+z_{12}\right| \ge 2, \\
\left| x_{11}-x_{12}+y_{11}-y_{12}+z_{11}-z_{12}\right| \ge 2, \\
\left| x_{11}-x_{14}+y_{11}-y_{14}-z_{11}+z_{14}\right| \ge 2, \\
\left| x_{12}-x_{13}+y_{12}-y_{13}-z_{12}+z_{13}\right| \ge 2, \\
\left| x_{12}-x_{14}+y_{12}-y_{14}+z_{12}-z_{14}\right| \ge 2, \\
\left| x_{13}-x_{14}+y_{13}-y_{14}+z_{13}-z_{14}\right| \ge 2.
\end{array}\right. \eqno (2.9)$$

\medskip
By optimizing $f_1(X)$ defined by $(2.7)$ under constraints $(2.8)$ and $(2.9)$ $($by Matlab$)$, one can deduce the following result.

\medskip\noindent
{\bf Theorem 2.1.} {\it If $O+X$ is a local packing with Figure $2.1$ as its color graph, then we have
$$f_1(X)\ge {38}/{27},$$
where the equality holds if $X$ is the set defined in Example $2.1.$ Consequently, the local density ${\rm vol}(O)/\Pi(X)$ attains $18/19$ as a local minimum.}

\medskip
\noindent
{\bf Experiment 2.2.} Assume that $C+X$ is a general local cuboctahedral packing associated with the color graph $G$ in Figure 2.2. Then $\Pi (X)$ has $14$ faces $Q_i$ determined by (2.1)
and $24$ vertices ${\bf p}_{i,j,k}=(x_{i,j,k}, y_{i,j,k}, z_{i,j,k})$ determined by (2.2), where $\{ i, j, k\}$ runs over the $24$ triples $\{1, 2, 4\}$, $\{1, 4, 5\}$, $\{1, 5, 14\}$,
$\{1, 2, 10\}$, $\{1, 10, 11\}$, $\{1, 11, 14\}$, $\{2, 3, 4\}$, $\{2, 3, 10\}$, $\{3, 4, 7\}$, $\{3, 7, 8\}$, $\{3, 8, 9\}$, $\{3, 9, 10\}$, $\{4, 5, 6\}$, $\{4, 6, 7\}$, $\{5, 6, 14\}$,
$\{6, 7, 8\}$, $\{6, 8, 13\}$, $\{6, 13, 14\}$, $\{8, 9, 12\}$, $\{8, 12, 13\}$, $\{9, 10,$ $12\}$, $\{10, 11, 12\}$, $\{11, 12, 14\}$ and $\{12, 13, 14\}$. Then, we rewrite the vertices of $Q_i$ as
$$\begin{array}{lllll}
Q_1:& {\bf p}_1^1={\bf p}_{1,2,10},&  {\bf p}_1^2={\bf p}_{1,10,11},&  {\bf p}_1^3={\bf p}_{1,11,14},&  {\bf p}_1^4={\bf p}_{1,5,14},\\
\vspace{-0.36cm}
&&&&\\
& {\bf p}_1^5={\bf p}_{1,4,5},& {\bf p}_1^6={\bf p}_{1,2,4};& &\\
\vspace{-0.36cm}
&&&&\\
Q_2:& {\bf p}_2^1={\bf p}_{1,2,4},&  {\bf p}_2^2={\bf p}_{2,3,4},&  {\bf p}_2^3={\bf p}_{2,3,10},&  {\bf p}_2^4={\bf p}_{1,2,10};\\
\vspace{-0.36cm}
&&&&\\
Q_3:& {\bf p}_3^1={\bf p}_{2,3,4},&  {\bf p}_3^2={\bf p}_{3,4,7},&  {\bf p}_3^3={\bf p}_{3,7,8},&  {\bf p}_3^4={\bf p}_{3,8,9},\\
\vspace{-0.36cm}
&&&&\\
& {\bf p}_3^5={\bf p}_{3,9,10},& {\bf p}_3^6={\bf p}_{2,3,10};&&\\
\vspace{-0.36cm}
&&&&\\
Q_4:& {\bf p}_4^1={\bf p}_{1,2,4},&  {\bf p}_4^2={\bf p}_{1,4,5},&  {\bf p}_4^3={\bf p}_{4,5,6},&  {\bf p}_4^4={\bf p}_{4,6,7},\\
\vspace{-0.36cm}
&&&&\\
& {\bf p}_4^5={\bf p}_{3,4,7},& {\bf p}_4^6={\bf p}_{2,3,4};&&\\
\vspace{-0.36cm}
&&&&\\
Q_5:& {\bf p}_5^1={\bf p}_{1,4,5},&  {\bf p}_5^2={\bf p}_{1,5,14},&  {\bf p}_5^3={\bf p}_{5,6,14},&  {\bf p}_5^4={\bf p}_{4,5,6};\\
\vspace{-0.36cm}
&&&&\\
Q_6:& {\bf p}_6^1={\bf p}_{4,5,6},&  {\bf p}_6^2={\bf p}_{5,6,14},&  {\bf p}_6^3={\bf p}_{6,13,14},&  {\bf p}_6^4={\bf p}_{6,8,13},\\
\vspace{-0.36cm}
&&&&\\
& {\bf p}_6^5={\bf p}_{6,7,8},& {\bf p}_6^6={\bf p}_{4,6,7};&&\\
\vspace{-0.36cm}
&&&&\\
Q_7:& {\bf p}_7^1={\bf p}_{3,4,7},&  {\bf p}_7^2={\bf p}_{4,6,7},&  {\bf p}_7^3={\bf p}_{6,7,8},&  {\bf p}_7^4={\bf p}_{3,7,8};\\
\vspace{-0.36cm}
&&&&\\
Q_8:& {\bf p}_8^1={\bf p}_{3,7,8},&  {\bf p}_8^2={\bf p}_{6,7,8},&  {\bf p}_8^3={\bf p}_{6,8,13},&  {\bf p}_8^4={\bf p}_{8,12,13},\\
\vspace{-0.36cm}
&&&&\\
& {\bf p}_8^5={\bf p}_{8,9,12},& {\bf p}_8^6={\bf p}_{3,8,9};&&\\
\vspace{-0.36cm}
&&&&\\
Q_9:& {\bf p}_9^1={\bf p}_{3,8,9},&  {\bf p}_9^2={\bf p}_{8,9,12},&  {\bf p}_9^3={\bf p}_{9,10,12},&  {\bf p}_9^4={\bf p}_{3,9,10};\\
\vspace{-0.36cm}
&&&&\\
Q_{10}:& {\bf p}_{10}^1={\bf p}_{1,2,10},&  {\bf p}_{10}^2={\bf p}_{2,3,10},&  {\bf p}_{10}^3={\bf p}_{3,9,10},&  {\bf p}_{10}^4={\bf p}_{9,10,12},\\
\vspace{-0.36cm}
&&&&\\
& {\bf p}_{10}^5={\bf p}_{10,11,12},& {\bf p}_{10}^6={\bf p}_{1,10,11};&&\\
\vspace{-0.36cm}
&&&&\\
Q_{11}:& {\bf p}_{11}^1={\bf p}_{1,10,11},&  {\bf p}_{11}^2={\bf p}_{10,11,12},&  {\bf p}_{11}^3={\bf p}_{11,12,14},&  {\bf p}_{11}^4={\bf p}_{1,11,14};\\
\vspace{-0.36cm}
&&&&\\
Q_{12}:& {\bf p}_{12}^1={\bf p}_{8,9,12},&  {\bf p}_{12}^2={\bf p}_{8,12,13},&  {\bf p}_{12}^3={\bf p}_{12,13,14},&  {\bf p}_{12}^4={\bf p}_{11,12,14},\\
\vspace{-0.36cm}
&&&&\\
& {\bf p}_{12}^5={\bf p}_{10,11,12},& {\bf p}_{12}^6={\bf p}_{9,10,12};&&\\
\vspace{-0.36cm}
&&&&\\
Q_{13}:& {\bf p}_{13}^1={\bf p}_{6,8,13},&  {\bf p}_{13}^2={\bf p}_{6,13,14},&  {\bf p}_{13}^3={\bf p}_{12,13,14},&  {\bf p}_{13}^4={\bf p}_{8,12,13};\\
\vspace{-0.36cm}
&&&&\\
Q_{14}:& {\bf p}_{14}^1={\bf p}_{1,5,14},&  {\bf p}_{14}^2={\bf p}_{1,11,14},&  {\bf p}_{14}^3={\bf p}_{11,12,14},&  {\bf p}_{14}^4={\bf p}_{12,13,14},\\
\vspace{-0.36cm}
&&&&\\
& {\bf p}_{14}^5={\bf p}_{6,13,14},& {\bf p}_{14}^6={\bf p}_{5,6,14}.&&
\end{array}$$

In this case, we have
$$n(1)=n(3)=n(4)=n(6)=n(8)=n(10)=n(12)=n(14)=6$$
and
$$n(2)=n(5)=n(7)=n(9)=n(11)=n(13)=4.$$
Then the volume of $\Pi (X)$ is
$$f_2(X)={1\over {12}}\sum_{i=1}^{14}\sum_{j=1}^{n(i)} \sqrt{\bigl(x_i^2+y_i^2+z_i^2\bigr)\Bigl( \bigl\| {\bf p}_i^1, {\bf p}_i^{j+1}\bigr\|^2\cdot \bigl\| {\bf p}_i^1, {\bf p}_i^{j+2}\bigr\|^2-\bigl\langle {\bf p}_i^{j+2}-{\bf p}_i^1,{\bf p}_i^{j+1}-{\bf p}_i^1\bigr\rangle^2 \Bigr) }.\eqno (2.10)$$
Furthermore, the colors of the vertices of $G$ imply

$$\left\{\begin{array}{l}
\left| x_1\right| \ge 2, \\
\left| x_2+y_2+z_2\right| \ge 4, \\
\left| x_3-y_3-z_3\right| \ge 4, \\
\left| z_4\right| \ge 2, \\
\left| x_5-y_5+z_5\right| \ge 4, \\
\left| y_6\right| \ge 2, \\
\left| x_7+y_7-z_7\right| \ge 4, \\
\left| x_8\right| \ge 2, \\
\left| x_9-y_9+z_9\right| \ge 4, \\
\left| y_{10}\right| \ge 2,\\
\left| x_{11}+y_{11}-z_{11}\right| \ge 4,\\
\left| z_{12}\right| \ge 2,\\
\left| x_{13}+y_{13}+z_{13}\right| \ge 4,\\
\left| x_{14}-y_{14}-z_{14}\right| \ge 4,
\end{array}
\right.\eqno(2.11)$$
and the colors of the edges of $G$ imply
$$\left\{\begin{array}{l}
\left| x_1-x_2-y_1+y_2-z_1+z_2\right| \ge 4, \\
\left| x_1-x_4+y_1-y_4-z_1+z_4\right| \ge 4,  \\

\left| y_1-y_5\right| \ge 2, \\

\left| x_1-x_{10}-y_1+y_{10}+z_1-z_{10}\right| \ge 4, \\

\left| z_1-z_{11}\right| \ge 2,  \\

\left| x_1-x_{14}+y_1-y_{14}+z_1-z_{14}\right| \ge 4,  \\

\left| x_2-x_3\right| \ge 2, \\

\left| y_2-y_4\right| \ge 2, \\

\left| z_2-z_{10}\right| \ge 2, \\

\left| x_3-x_4-y_3+y_4+z_3-z_4\right| \ge 4, \\

\left| y_3-y_7\right| \ge 2, \\

\left| x_3-x_8+y_3-y_8+z_3-z_8\right| \ge 4, \\

\left| z_3-z_9\right| \ge 2, \\

\left| x_3-x_{10}+y_3-y_{10}-z_3+z_{10}\right| \ge 4, \\

\left| x_4-x_5-y_4+y_5-z_4+z_5\right| \ge 4, \\

\left| x_4-x_6+y_4-y_6+z_4-z_6\right| \ge 4, \\

\left| x_4-x_7\right| \ge 2, \\

\left| x_5-x_6\right| \ge 2, \\

\left| z_5-z_{14}\right| \ge 2, \\

\left| x_6-x_7-y_6+y_7-z_6+z_7\right| \ge 4, \\

\left| x_6-x_8-y_6+y_8+z_6-z_8\right| \ge 4, \\

\left| z_6-z_{13}\right| \ge 2, \\

\left| x_6-x_{14}+y_6-y_{14}-z_6+z_{14}\right| \ge 4, \\

\left| z_7-z_8\right| \ge 2, \\

\left| y_8-y_9\right| \ge 2, \\

\left| x_8-x_{12}+y_8-y_{12}-z_8+z_{12}\right| \ge 4, \\

\left| x_8-x_{13}-y_8+y_{13}-z_8+z_{13}\right| \ge 4, \\

\left| x_9-x_{10}\right| \ge 2, \\

\left| x_9-x_{12}-y_9+y_{12}-z_9+z_{12}\right| \ge 4, \\

\left| x_{10}-x_{11}-y_{10}+y_{11}-z_{10}+z_{11}\right| \ge 4, \\

\left| x_{10}-x_{12}+y_{10}-y_{12}+z_{10}-z_{12}\right| \ge 4, \\

\left| x_{11}-x_{12}\right| \ge 2, \\

\left| y_{11}-y_{14}\right| \ge 2, \\

\left| y_{12}-y_{13}\right| \ge 2, \\

\left| x_{12}-x_{14}-y_{12}+y_{14}+z_{12}-z_{14}\right| \ge 4, \\

\left| x_{13}-x_{14}\right| \ge 2.
\end{array}\right. \eqno (2.12)$$

Similar to Theorem 2.1, by optimizing $f_2(X)$ defined by (2.10) under constrains (2.11) and (2.12), one can obtain the following result.

\medskip\noindent
{\bf Theorem 2.2.} {\it If $C+X$ is a local packing with Figure $2.2$ as its color graph, one can deduce that
$$f_2(X)\ge {{196}\over {27}},$$
where the equality holds if $X$ is the set defined in Example $2.2.$ Consequently, the local density ${\rm vol}(C)/\Pi(X)$ attains $45/49$ as a local minimum.}

\vspace{1cm}
\centerline{\Large\bf 3. Triangulated Color Graphs}

\bigskip
\noindent
It follows by Corollary 1.1 that, for any fixed centrally symmetric three-dimensional polytope $P$, to determine the minimum volume of the local cells $\Pi (X)$ or the maximum local density $\delta (P)$
for all possible packings $P+X$ it is sufficient to deal with the general local ones with ${\rm card}\{ X\} \le 26^3$ (as it was mentioned in Remark 1.2, this upper bound should be much improved).
Therefore, one can divide all those general local packings into finite classes according to their triangulated color graphs. Clearly, each class can be treated by computer as an optimization problem.
If the number of the classes is not too large, the whole problem can be treated by a computer.

\vspace{0.6cm}\noindent
{\bf 3.1. Triangulated Graphs}

\bigskip
\noindent
Let $n$, $e$ and $f$ denote the numbers of the vertices, the edges and the faces of a triangulated graph $G$, respectively. By {\it Euler's formula} we have
$$\left\{\begin{array}{l}
e=3(n-2),\\
f=2(n-2).
\end{array}\right.\eqno(3.1)$$

Let $g(n)$ denote the number of all the isomorphically distinct triangulated graphs (no color yet) with $n$ vertices, and let $\mathcal{G}_n$ denote a set of $g(n)$ isomorphically distinct
triangulated graphs with $n$ vertices. To determine the values of $g(n)$ and to generate a graph family $\mathcal{G}_n$ is a challenging job. In 1962, W. T. Tutte \cite{tutte} proved the following result.

\medskip\noindent
{\bf Lemma 3.1.}
$${{(4n-11)!}\over {6(n-2)\cdot (3n-7)!\cdot (n-2)!}}\le g(n)\le {{2\cdot (4n-11)!}\over {(3n-7)!\cdot (n-2)!}}.$$

\bigskip
By Stiring's formula, roughly speaking, this lemma means that
$${c\over {6 n^{3.5}}}\left({{256}\over {27}}\right)^n\ll g(n)\ll {c\over {n^{2.5}}}\left({{256}\over {27}}\right)^n,$$
where $c$ is a constant defined by
$$c={{3^{6.5}}\over { 2^{20.5}\pi^{0.5}}}\approx {{48}\over {10^4}}.$$

\medskip
Next, we will propose a process to generate $\mathcal{G}_n$ from $\mathcal{G}_{n-1}$. Our process is based on the following observation.

\medskip\noindent
{\bf Lemma 3.2.} {\it Assume that $G$ is a triangulated graph with $n$ vertices $(n\ge 4)$. If $G$ has exact $\varpi_k$ vertices of degree $k$, then we have
$\varpi_3\not= 0$, $\varpi_4\not= 0$ or $\varpi_5\not= 0$.}

\medskip\noindent
{\bf Proof.} First of all, it is easy to see that
$$\sum_{j=3}^{n-1}\varpi_j=n.\eqno(3.2)$$
On the other hand, by the first equation of (3.1) and double counting the edges, it can be deduced that
$$\sum_{j=3}^{n-1}j\cdot\varpi_j=6(n-2).\eqno(3.3)$$

If, on the contrary, $\varpi_3=\varpi_4=\varpi_5=0$, then (3.3) and (3.2) imply that
$$6(n-2)=\sum_{j=6}^{n-1}j\cdot\varpi_j \ge 6 \sum_{j=6}^{n-1}\varpi_j =6n,$$
which is apparently impossible. Lemma 3.2 is proved. \hfill{$\Box$}

\bigskip
Next, we observe three degenerating processes corresponding to the three inequalities of Lemma 3.2, respectively.

\medskip\noindent
{\bf Case 1.} $\varpi_3\not= 0$. If the degree of ${\bf v}_i$ is three in $G$. Assume that $E_i^1$, $E_i^2$ and $E_i^3$ are the three edges meeting at ${\bf v}_i$, and ${\bf v}_i^j$ is the other
end of $E_i^j$, then by deleting all ${\bf v}_i$, $E_i^1$, $E_i^2$ and $E_i^3$ (as illustrated by Figure 3.1) we get a triangulated graph with $n-1$ vertices.

\begin{figure}[ht]
\centering
\includegraphics[height=3.5cm,width=10cm,angle=0]{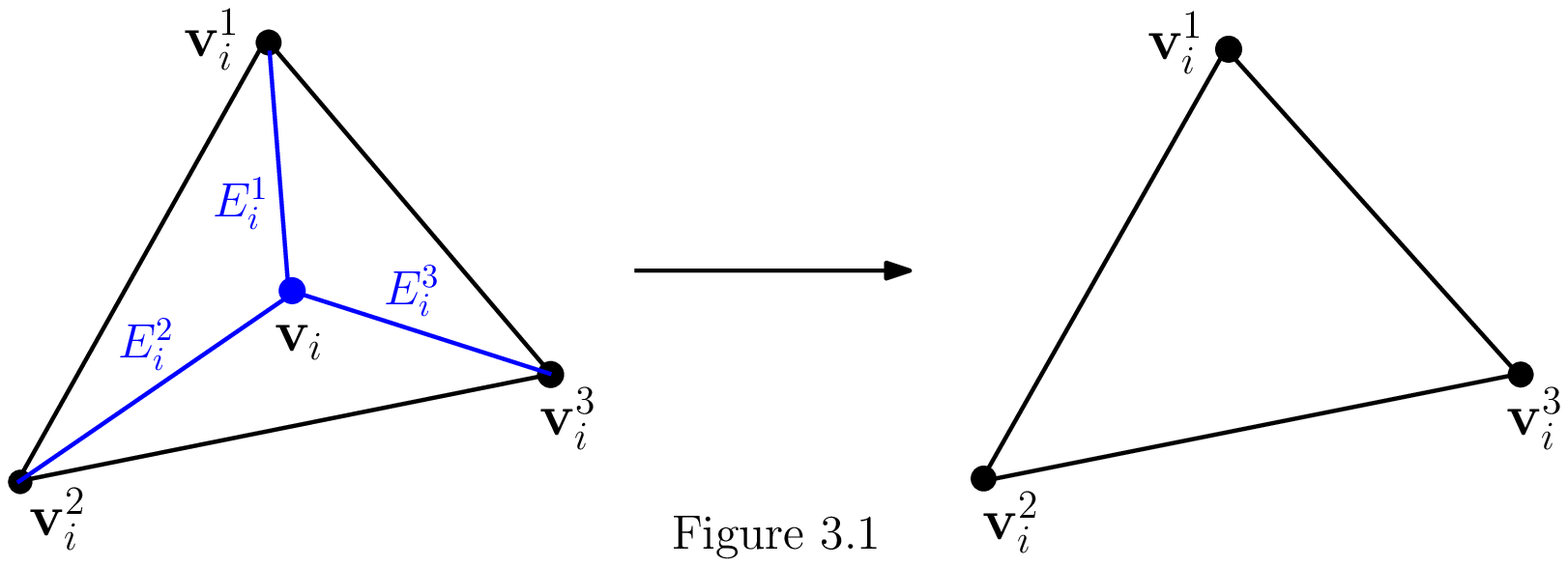}
\end{figure}

\noindent
{\bf Case 2.} $\varpi_4\not= 0$. If the degree of ${\bf v}_i$ is four in $G$. Assume that $E_i^1$, $E_i^2$, $E_i^3$ and $E_i^4$ (in an anti-clock order) are the four edges meeting at ${\bf v}_i$,
and ${\bf v}_i^j$ is the other end of $E_i^j$, then we delete all ${\bf v}_i$, $E_i^1$, $E_i^2$, $E_i^3$ and $E_i^4$, and add a new edge connecting ${\bf v}_i^2$ and ${\bf v}_i^4$
(the local change of the graph is illustrated by Figure 3.2). In this way, we get a new graph in $\mathcal{G}_{n-1}$.

\smallskip
\begin{figure}[ht]
\centering
\includegraphics[height=4cm,width=10.5cm,angle=0]{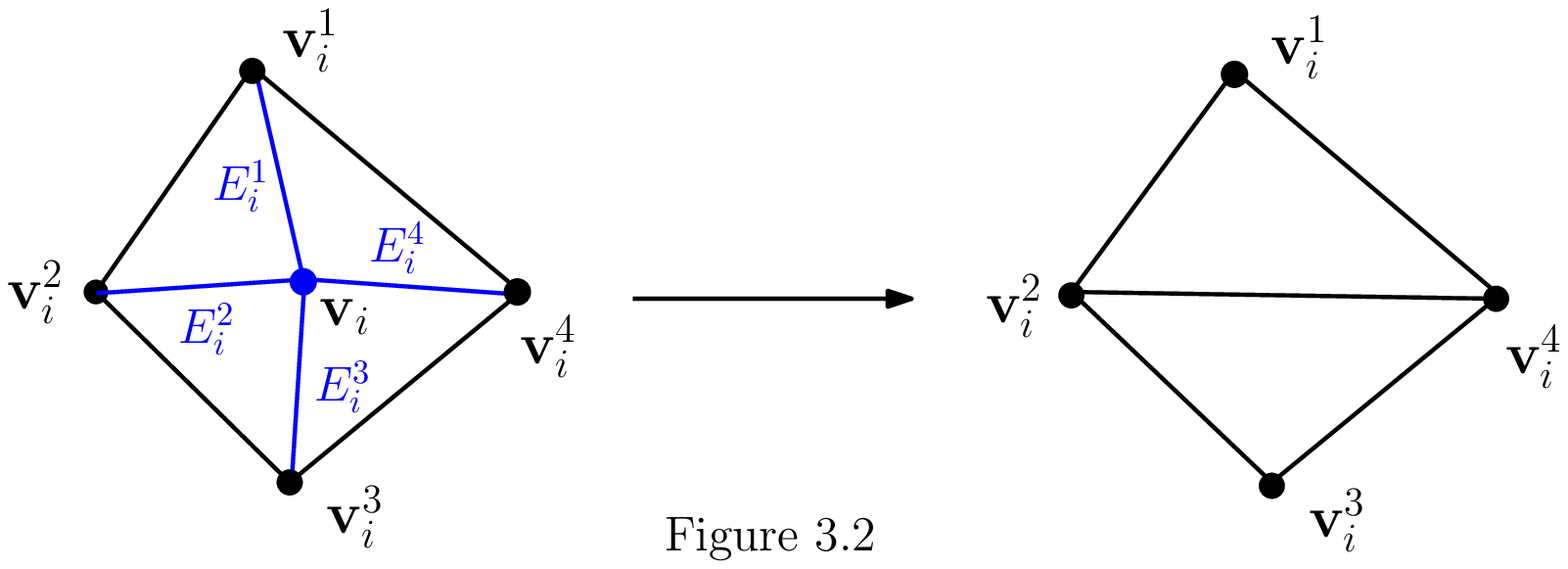}
\end{figure}

\medskip\noindent
{\bf Case 3.} $\varpi_5\not= 0$. If the degree of ${\bf v}_i$ is five in $G$. Assume that $E_i^1$, $E_i^2$, $E_i^3$, $E_i^4$ and $E_i^5$ (in an anti-clock order) are the five edges meeting at ${\bf v}_i$, and ${\bf v}_i^j$ is the other end of $E_i^j$, then we delete all ${\bf v}_i$, $E_i^1$, $E_i^2$, $E_i^3$, $E_i^4$ and $E_i^5$ and add two new edges connecting ${\bf v}_i^1$ with ${\bf v}_i^3$ and
${\bf v}_i^4$, respectively (the local change of the graph is illustrated by Figure 3.3). In this way, we get a new graph in $\mathcal{G}_{n-1}$.

\medskip
\begin{figure}[ht]
\centering
\includegraphics[height=4cm,width=11.3cm,angle=0]{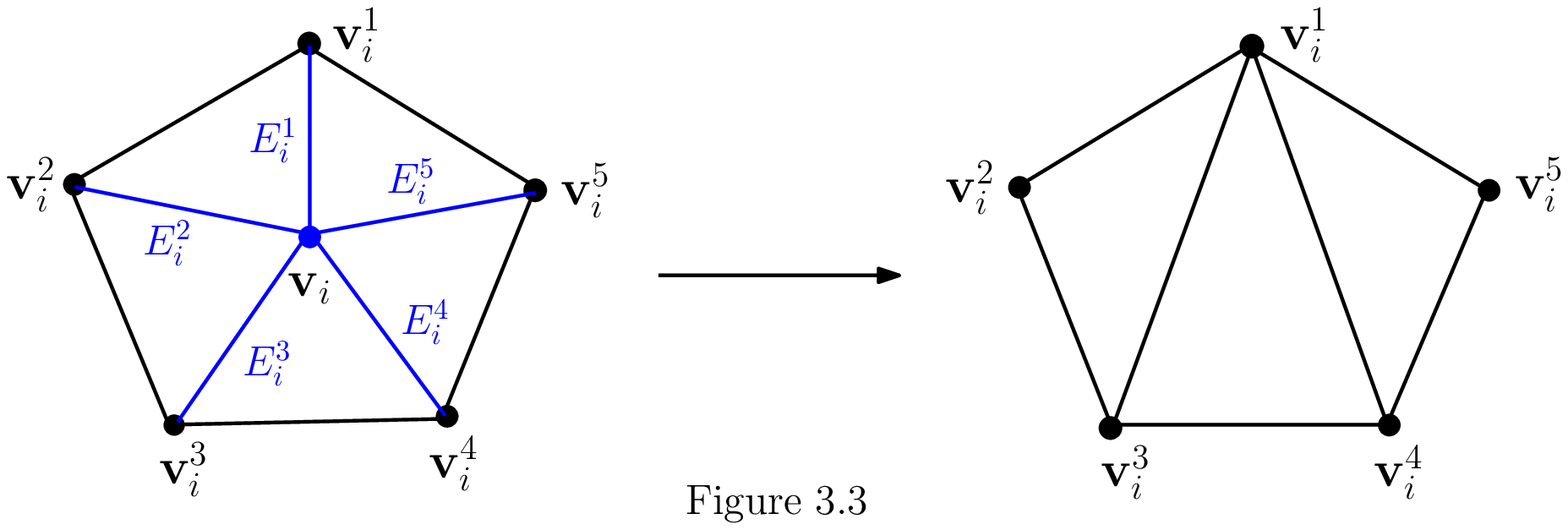}
\end{figure}

Based on this observation, we propose the following process to generate $\mathcal{G}_n$ from $\mathcal{G}_{n-1}$. In fact, what we will process is just the inverse of the previous
three cases and what we will get is a set $\mathcal{G}^*_n$ containing $\mathcal{G}_n$.

Assume that $\mathcal{G}_{n-1}=\{ G_1, G_2, \ldots, G_\kappa\}$, where $G_i$ are triangulated graphs with $n-1$ vertices, $3(n-3)$ edges, and $2(n-3)$ triangular faces, and $\kappa =g(n-1)$.

Taking a graph $G_i\in \mathcal{G}_{n-1}$, for convenience, we enumerate its faces by $\triangle_i^1$, $\triangle_i^2$, $\ldots $, $\triangle_i^{2(n-3)}$, and enumerate its edges by $E_i^1$, $E_i^2$,
$\ldots $, $E_i^{3(n-3)}$.

\smallskip\noindent
{\bf Generating process 1.} Adding a new vertex ${\bf v}_n$ in the interior of $\triangle_i^j$ and connecting it with the three vertices of $\triangle_i^j$ by three new edges, we get a triangulated graph
$G_i^{1,j}$. In this way, we obtain $2\kappa \cdot (n-3)$ triangulated graphs with $n$ vertices. For convenience, let $\mathcal{G}_n^1$ denote the set of these graphs.

\smallskip\noindent
{\bf Generating process 2.} Assume that two triangular faces $\triangle$ and $\triangle '$ joining at the entire edge $E_i^j$, ${\bf v}$ is the only vertex of $\triangle$ which is not an end of $E_i^j$ and
${\bf v}'$ is the only vertex of $\triangle '$ which is not an end of $E_i^j$. Adding a new vertex ${\bf v}_n$ in the relative interior of $E_i^j$ and connecting it with ${\bf v}$ and ${\bf v}'$ by two new edges, we get a triangulated graph $G_i^{2,j}$. In this way, we obtain $3\kappa \cdot (n-3)$ triangulated graphs with $n$ vertices. For convenience, let $\mathcal{G}_n^2$ denote the set of these graphs.

\smallskip\noindent
{\bf Generating process 3.} Assume that two triangular faces $\triangle$ and $\triangle '$ joining $\triangle_i^j$ at entire edges, respectively, ${\bf v}$ is the vertex of $\triangle$ which is not a vertex of $\triangle_i^j$ and
${\bf v}'$ is the vertex of $\triangle '$ which is not a vertex of $\triangle_i^j$. Deleting the two edges $\triangle\cap \triangle_i^j$ and $\triangle '\cap \triangle_i^j$, adding a new vertex ${\bf v}_n$ in the interior region of $\triangle_i^j$ and connecting it with the five vertices of $\triangle\cup \triangle_i^j\cup \triangle'$ by five new edges, we get a triangulated graph $G_i^{3,j}$. In this way, we obtain $6\kappa \cdot (n-3)$ triangulated graphs with $n$ vertices. For convenience, let $\mathcal{G}_n^3$ denote the set of these graphs.

\medskip
As a consequence of Lemma 3.2, we obtained the following result.

\smallskip
\noindent
{\bf Lemma 3.3.} {\it For $\mathcal{G}_n$ and $\mathcal{G}_{n-1}$, we have
$$\mathcal{G}_n\subseteq \mathcal{G}_n^1\cup\mathcal{G}_n^2\cup \mathcal{G}_n^3$$
and}
$${{g(n)}\over {g(n-1)}}={{{\rm card}\{ \mathcal{G}_n\}}\over {{\rm card}\{ \mathcal{G}_{n-1}\}}} \le 11(n-3).$$

\vspace{0.6cm}\noindent
{\bf 3.2. Triangulated Color Graphs}

\medskip
\noindent
Assume that $P$ has $m$ pairs of faces and let us consider the general local packings $P+X$ with ${\rm card}(X)=n+1$. It follows from (3.1) that the triangulated color graph $G$ corresponding to $P+X$
has $n$ vertices, $3(n-2)$ edges and $2(n-2)$ triangular faces. Furthermore, the vertices and the edges of $G$ are colored by $m$ different colors.

Let $\varphi (n)$ denote the number of the isomorphically distinct triangulated color graphs corresponding to general local packings $P+X$ with ${\rm card}(X)=n+1$ and let $\mathcal{G}$ denote a set
of all possible isomorphically distinct triangulated color graphs corresponding to some general
local packings $P+X$ satisfying ${\rm vol}(\Pi (X))= \omega (P)$. By Lemma 3.1 and Corollary 1.1 one can easily deduce the following upper bounds for $\varphi (n)$ and ${\rm card}\{ \mathcal{G}\}.$

\medskip\noindent
{\bf Lemma 3.4.} {\it For a general centrally symmetric polytope $P\in \mathcal{K}$ with $m$ pairs of faces we have
$$\varphi (n)\le g(n)\cdot m^n\cdot m^{3(n-2)}\le {{2\cdot (4n-11)!\cdot m^{4n-6}}\over {(3n-7)!\cdot (n-2)!}}$$
and}
$${\rm card}\{ \mathcal{G}\}\le \sum_{n=4}^{26^3}\varphi (n)\le \sum_{n=4}^{26^3} {{2\cdot (4n-11)!\cdot m^{4n-6}}\over {(3n-7)!\cdot (n-2)!}}.$$

\smallskip\noindent
{\bf Corollary 3.1.} When $P=O$, we have $m=4$, $\tau (O)\le 10$ (Remark 1.2) and
$${\rm card}\{\mathcal{G} \} \le \sum_{n=4}^{11^3} {{2\cdot (4n-11)!\cdot 4^{4n-6}}\over {(3n-7)!\cdot (n-2)!}}.$$
When $P=C$, we have $m=7$, $\tau (C)\le 10$ (Remark 1.2) and
$${\rm card}\{\mathcal{G} \} \le \sum_{n=4}^{11^3} {{2\cdot (4n-11)!\cdot 7^{4n-6}}\over {(3n-7)!\cdot (n-2)!}}.$$

\vspace{0.6cm}\noindent
{\bf 3.3. The Adjacency Matrix of a Triangulated Color Graph}

\medskip
\noindent
For a computer, matrices should be more recognizable than color graphs. Therefore, in this section we will introduce an {\it adjacency matrix} for a color graph.

\smallskip\noindent
{\bf Definition 3.1.} Assume that $G$ is a triangulated color graph with $n$ vertices ${\bf v}_1$, ${\bf v}_2$, $\ldots,$ ${\bf v}_n$ and $m$ colors
${\bf w}_1$, ${\bf w}_2$, $\ldots ,$ ${\bf w}_m$. We define its $n\times n$ adjacency matrix $M=\left(a_{i,j}\right)$ by
$$a_{ij}=\left\{
\begin{array}{ll}
s & \mbox{if $i=j$ and the vertex ${\bf v}_i$ is in color ${\bf w}_s$;}\\
t & \mbox{if $i\not=j$ and there is a ${\bf w}_t$ color edge connecting ${\bf v}_i$ and ${\bf v}_j$;}\\
0 & \mbox{if $i\not=j$, but ${\bf v}_i$ and ${\bf v}_j$ are not connected.}
\end{array}\right.$$

\smallskip\noindent
{\bf Example 3.1.} The adjacency matrix of the triangulated color graph of $O+X$ illustrated by Figure 2.1 is

$$\left(\begin{array}{cccccccccccccc}
1&2&0&3&4&0&0&0&0&1&3&0&0&4\\
2&4&1&4&0&0&0&0&0&3&0&0&0&0\\
0&1&2&1&0&0&4&4&3&3&0&0&0&0\\
3&4&1&3&2&4&1&0&0&0&0&0&0&0\\
4&0&0&2&1&1&0&0&0&0&0&0&0&3\\
0&0&0&4&1&4&2&1&0&0&0&0&3&3\\
0&0&4&1&0&2&3&3&0&0&0&0&0&0\\
0&0&4&0&0&1&3&1&4&0&0&3&2&0\\
0&0&3&0&0&0&0&4&1&1&0&2&0&0\\
1&3&3&0&0&0&0&0&1&4&2&4&0&0\\
3&0&0&0&0&0&0&0&0&2&3&1&0&4\\
0&0&0&0&0&0&0&3&2&4&1&3&4&1\\
0&0&0&0&0&3&0&2&0&0&0&4&4&1\\
4&0&0&0&3&3&0&0&0&0&4&1&1&2
\end{array}\right).
$$
In this case, we have $n=14$ and $m=4$, where ${\bf w}_1$ is red, ${\bf w}_2$ is green, ${\bf w}_3$ is blue, and ${\bf w}_4$ is black.

\smallskip\noindent
{\bf Remark 3.1.} The triangular property of the graph can be traced in the adjacency matrix, but much less obvious.

\smallskip\noindent
{\bf Example 3.2.} The adjacency matrix of the triangulated color graph of $C+X$ illustrated by Figure 2.2 is

$$\left(\begin{array}{cccccccccccccc}
1&6&0&7&2&0&0&0&0&5&3&0&0&4\\
6&4&1&2&0&0&0&0&0&3&0&0&0&0\\
0&1&6&5&0&0&2&4&3&7&0&0&0&0\\
7&2&5&3&6&4&1&0&0&0&0&0&0&0\\
2&0&0&6&5&1&0&0&0&0&0&0&0&3\\
0&0&0&4&1&2&6&5&0&0&0&0&3&7\\
0&0&2&1&0&6&7&3&0&0&0&0&0&0\\
0&0&4&0&0&5&3&1&2&0&0&7&6&0\\
0&0&3&0&0&0&0&2&5&1&0&6&0&0\\
5&3&7&0&0&0&0&0&1&2&6&4&0&0\\
3&0&0&0&0&0&0&0&0&6&7&1&0&2\\
0&0&0&0&0&0&0&7&6&4&1&3&2&5\\
0&0&0&0&0&3&0&6&0&0&0&2&4&1\\
4&0&0&0&3&7&0&0&0&0&2&5&1&6
\end{array}\right).
$$
In this case, we have $n=14$ and $m=7$, where ${\bf w}_1$ is red, ${\bf w}_2$ is green, ${\bf w}_3$ is blue, ${\bf w}_4$ is yellow, ${\bf w}_5$ is purple, ${\bf w}_6$ is brown, and ${\bf w}_7$ is black.

\vspace{1cm}
\centerline{\Large\bf 4. Localizations with Truncaters}

\vspace{0.6cm}\noindent
In principle, Conjecture Z can be proved or disproved by dealing with a big number of optimizations of (2.6) type based on Corollary 1.1 and Corollary 3.1. However, the number of the cases is too big, even for the computer. For the purpose to reduce the number of cases, we introduce a truncater next.

\vspace{0.6cm}\noindent
{\bf 4.1. Truncaters}

\medskip
\noindent
We recall that $K$ is a three-dimensional centrally symmetric convex body centered at the origin satisfying (1.1), $\mathfrak{X}$ is the family of all discrete sets $X$ such that ${\bf o}\in X$ and $K+X$ is a packing. Let $\Gamma$ be a {\it truncater} (a centrally symmetric convex body centered at the origin). We define
$$\omega (\Gamma, K)=\min_{X\in \mathfrak{X}}{\rm vol}(\Gamma \cap \Pi (X)).$$
It is obvious that
$$\omega (\Gamma, K) \le \omega (K)$$
holds for all truncaters $\Gamma $. Therefore, for a particular $K$ if we can luckily choose a suitable truncater $\Gamma $ such that, for any polytope $P$, the volume ${\rm vol}(\Gamma \cap P)$ is relatively easy to compute and
$$\omega (\Gamma, K) = \omega (K),$$
we will be able to determine the value of $\delta (K)$. Furthermore, if there is a lattice packing $K+\Lambda$ attending $\delta (K)$, then one can determine the values of $\delta^t(K)$ and $\delta^l(K)$ by
$$\delta^t(K)=\delta^l(K)=\delta (K).$$

Assume that $K+X$ is a local packing attending the maximal local density $\delta (K)$, where $X=\{ {\bf o}, {\bf x}_1, {\bf x}_2, \ldots , {\bf x}_n\}$. It is easy to see that
$${\rm card}\{ X\}\ge 5.$$
Otherwise, all ${\bf x}_i$ will lie on one side of (or on) certain hyperplane passing the origin. Consequently, the Dirichlet-Voronoi cell $\Pi (X)$ would be unbounded and therefore the packing can not attend the minimal $\delta (K)$.
We recall that $H_{\bf x}$ is the bisector of ${\bf o}$ and ${\bf x}$. If $\Gamma \cap H_{{\bf x}_i}=\emptyset $, then we have
$$\Gamma \cap \Pi (X)=\Gamma \cap \Pi (X\setminus \{ {\bf x}_i\}).$$
For this reason, we introduce the following notion.

\medskip\noindent
{\bf Definition 4.1.} For a centrally symmetric convex body $\Gamma $ centered at the origin we define
$$\Gamma^* =\{ {\bf x}: \ \Gamma \cap H_{\bf x}\not = \emptyset \}$$
to be its colony.

\medskip
It is easy to show that, for every centrally symmetric convex body $\Gamma $ we have
$$2\Gamma \subseteq \Gamma^*\eqno (4.1)$$
and
$$\Gamma_1^*\subseteq \Gamma_2^*\eqno (4.2)$$
whenever $\Gamma_1\subseteq \Gamma_2$. In particular, if ${\bf x}\not\in \Gamma^*$, the translate $K+{\bf x}$ has no effect on $\Gamma \cap \Pi (X)$. Therefore, it is both interesting and useful to have a close look at $\Gamma^*$.

\medskip\noindent
{\bf Example 4.1.} If $\Gamma =\{ (x,y): \ |x| \le k, \ |y|\le 1\}$ and let $\gamma (x,y)$ denote the maximum number $\gamma $ such that $(\gamma x,\gamma y)\in \Gamma^*$, then by a routine computation one can deduce that
$$\gamma (x,y)=\left\{
\begin{array}{ll}
{{2(k|x|+1)}\over {x^2+1}}& \mbox{if $|y|=1$,}\\

\vspace{-0.35cm}
&\\
{{2(k^2+|y|)}\over {k^2+y^2}}& \mbox{if $|x|=k$.}\\
\end{array}
\right.$$
In particular, when $k=1$, the colony of a square is illustrated by Figure 4.1, which is no longer a convex domain.

\begin{figure}[ht]
\centering
\includegraphics[height=5.5cm,width=5.5cm,angle=0]{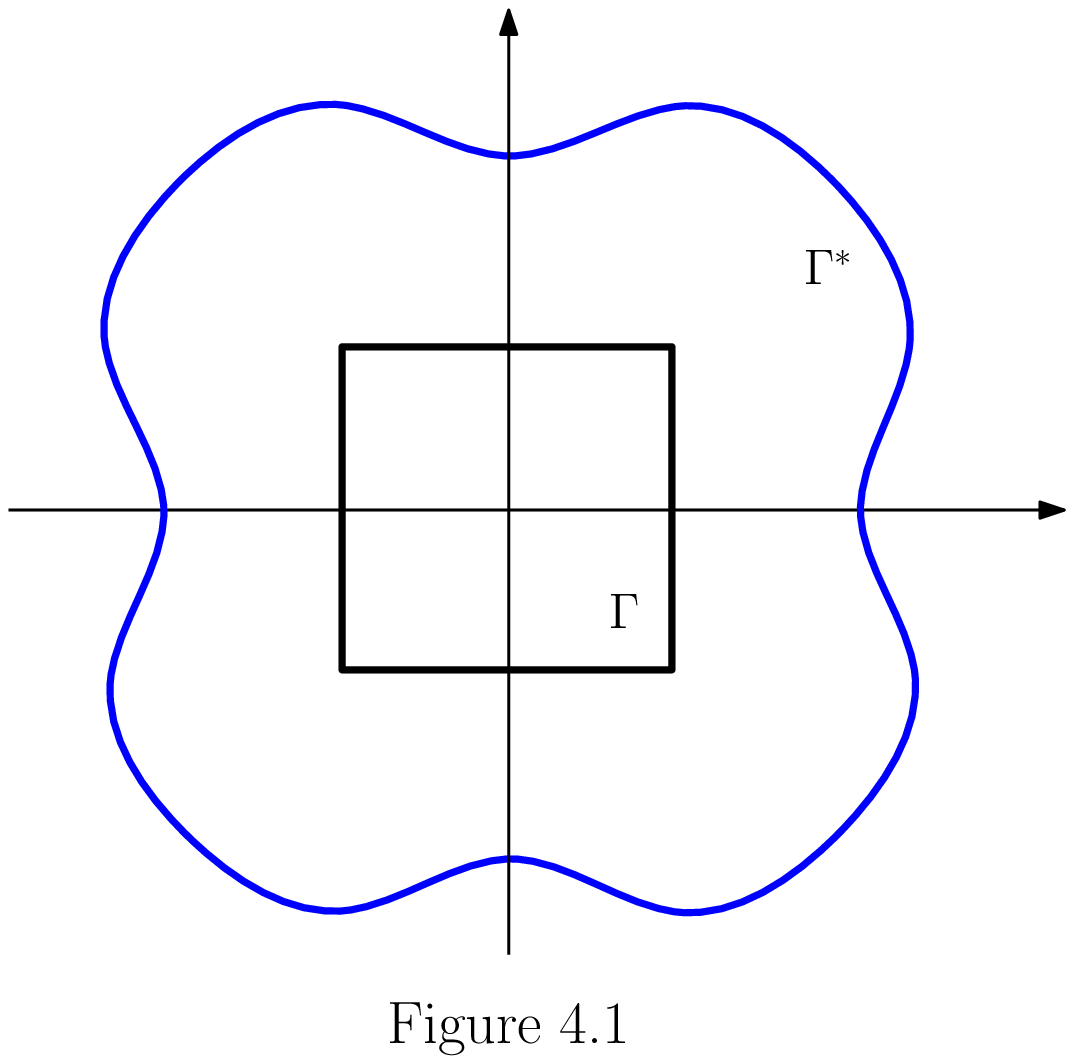}
\end{figure}

To characterize the equality case in (4.1) we have the following result.

\medskip\noindent
{\bf Theorem 4.1.} {\it An $n$-dimensional centrally symmetric convex body $\Gamma $ satisfies $\Gamma^*=2\Gamma $ if and only if $\Gamma $ is a ball centered at the origin.}

\medskip\noindent
{\bf Proof.} The if part is trivial. Now we deal with the only if part.

First, let us reduce the problem to two dimensions. Let $H$ be a two-dimensional hyperplane passing the origin. On one hand, by (4.1) we have
$$2\Gamma \cap H\subseteq (\Gamma \cap H)^*.$$
On the other hand, by (4.2) and $\Gamma^*=2\Gamma $ we have
$$(\Gamma \cap H)^*\subseteq \Gamma^*\cap H=2\Gamma \cap H.$$
As a conclusion, if $\Gamma^*=2\Gamma,$ for every two-dimensional hyperplane $H$ passing the origin we have
$$(\Gamma \cap H)^*=2\Gamma \cap H.$$
Therefore, to prove the theorem, it is sufficient to show the two-dimensional case only. That is, we may assume that $\Gamma $ is a two-dimensional centrally symmetric convex domain centered at the origin satisfying $\Gamma^*=2\Gamma $.

Second, we claim that $\Gamma $ is differentiable at any point ${\bf x}=(x,y)\in \partial (\Gamma ).$ In other words, $\Gamma $ has unique tangent line at every boundary point.

\begin{figure}[ht]
\centering
\includegraphics[height=4.5cm,width=6cm,angle=0]{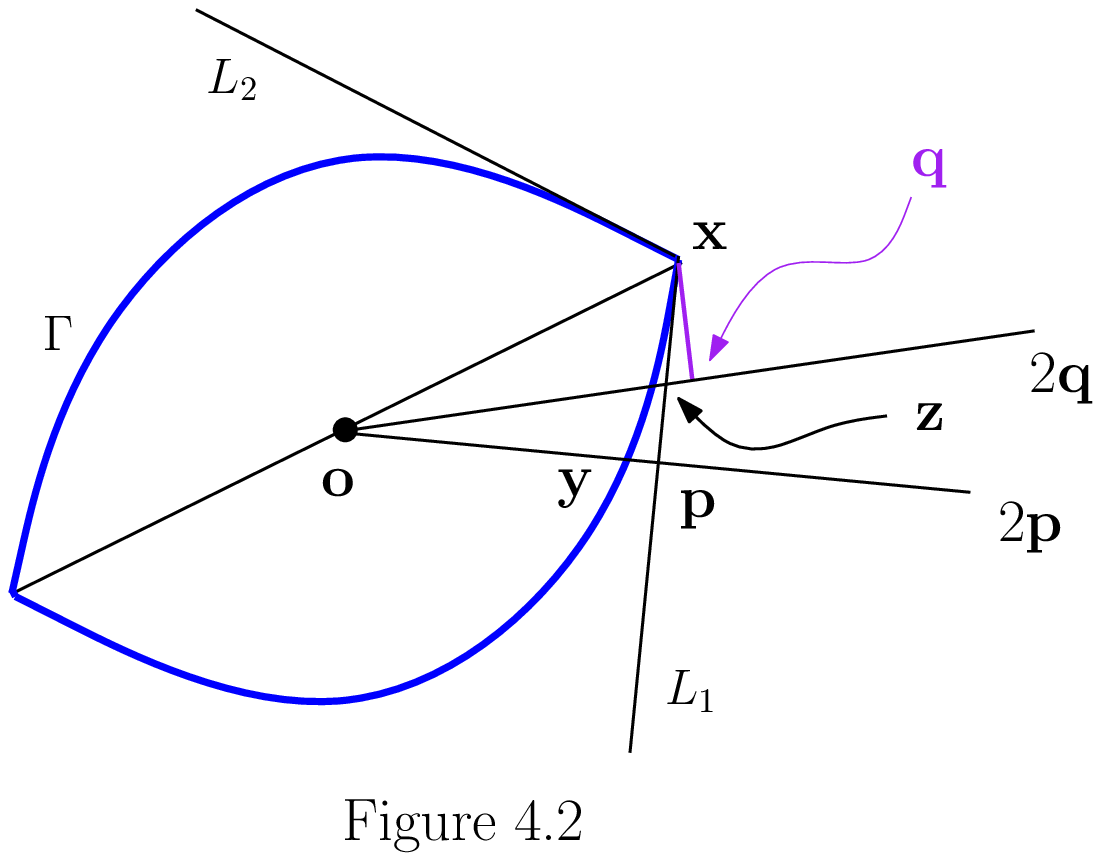}
\end{figure}

If, on the contrary, $\Gamma $ has two different tangent lines $L_1$ and $L_2$ at ${\bf x}\in \partial (\Gamma )$, as shown by Figure 4.2. We may assume that the angle between ${\bf xo}$ and $L_1$
is less than $\pi /2$. Therefore, there is a point ${\bf p}\in L_1$ such that ${\bf op}$ is perpendicular to $L_1$ and $\partial (\Gamma )$ has a unique point ${\bf y}$ on the segment ${\bf op}$.
By the definition of $\Gamma^*$ it is clear that $2{\bf p}\in \partial (\Gamma^*)$. Then it follows that ${\bf y}$ and ${\bf p}$ are identical and thus the whole segment ${\bf xp}$ belongs to $\partial (\Gamma )$.
Take ${\bf z}$ to be the middle point of ${\bf xp}$ and choose ${\bf q}$ to be the point on the straight line $L$ passing both ${\bf o}$ and ${\bf z}$ such that ${\bf xq}$ is perpendicular to $L$. Clearly,
${\bf z}$ and ${\bf q}$ are not identical and thus $2{\bf z}\not= 2{\bf q}$. However, we have $2{\bf q}\in \partial (\Gamma^*)$. Thus, we have obtained
$$\Gamma^*\not= 2\Gamma ,$$
which contradicts the assumption.

The rest of the proof is similar to the proof of Theorem 1.1. Assume that ${\bf x}=(x,y)\in \partial (\Gamma )$ is defined by
$$\left\{
\begin{array}{ll}
x=r(\theta )\cos \theta ,&\\
y=r(\theta )\sin \theta ,&
\end{array}\right.$$
where $r(\theta ) >0$ and $0\le \theta < 2\pi $. It can be deduced that $r(\theta )$ is differentiable. In fact, if the tangent ratio of $\partial (\Gamma)$ at ${\bf x}$ is $k(\theta)$, one can deduce that
$$r'(\theta )={{r(\theta ) \bigl(k(\theta )\sin \theta +\cos \theta \bigr)}\over {k(\theta )\cos \theta -\sin\theta }}.$$
Then the tangent direction ${\bf v}=(x', y')$ of $\Gamma $ at ${\bf x}$ is defined by
$$\left\{
\begin{array}{ll}
x'=r'(\theta )\cos \theta -r(\theta )\sin \theta,&\\
y'=r'(\theta )\sin \theta +r(\theta )\cos \theta.&
\end{array}\right.$$
Since ${\bf x}$ as a vector is an norm of $\Gamma $ at ${\bf x}$, we have $\langle {\bf x}, {\bf v}\rangle =0$ and therefore
$$r(\theta )\cos \theta \left(r'(\theta )\cos \theta -r(\theta )\sin \theta\right) +r(\theta )\sin \theta \left(r'(\theta )\sin \theta +r(\theta )\cos \theta\right)=0,$$

$$r(\theta )r'(\theta )\left(\sin^2\theta +\cos^2\theta \right)=r(\theta )r'(\theta )=0$$
and finally
$$r'(\theta )=0.$$
This means that $r(\theta )$ is a constant and therefore $\Gamma$ must be a circular domain. The theorem is proved. \hfill{$\square$}

\medskip\noindent
{\bf Example 4.2.} Assume that $\Lambda_1$ is the lattice defined in Example 1.1, it is easy to check that all the distances between the vertices and the center of $\Pi(\Lambda_1)$ are
$$r_1={{\sqrt{2033}}\over {57}}\approx 0.79103\ldots .$$ Therefore, to determine $\delta (O)$ it is reasonable to take $r_1B$ as a truncater.

\medskip\noindent
{\bf Example 4.3.} Assume that $\Lambda_2$ is the lattice defined in Example 1.2, it is easy to check that all the distances between the vertices and the center of $\Pi (\Lambda_2)$ are
$$r_2={{\sqrt{830}}\over {21}}\approx 1.37189\ldots .$$ Therefore, to determine $\delta (C)$ it is reasonable to take $r_2B$ as a truncater.

\vspace{0.6cm}\noindent
{\bf 4.2. The Effective Neighbours}

\medskip
\noindent
To determine or estimate the value of $\omega (\Gamma, K)$, it is sufficient to study the local packings $K+X$ satisfying
$$X\subset \Gamma^*.\eqno (4.3)$$
For fixed $\Gamma$ and $K$, we call $K+{\bf x}$ an {\it effective neighbour} of $K$ with respect to the truncater $\Gamma$ if ${\bf x}\in \Gamma^*$ and ${\rm int}(K)\cap (K+{\bf x})=\emptyset $.
Let $m(\Gamma, K)$ to be the maximal number of effective neighbours of $K$ with respect to $\Gamma$ in a local packing $K+X$. In other words, to study $\omega (\Gamma, K)$ it is sufficient to study the local
packings with at most $m(\Gamma, K)$ neighbours.

\medskip\noindent
{\bf Theorem 4.2.} {\it Assume that both $\Gamma$ and $K$ are $n$-dimensional centrally symmetric convex bodies, then}
$$m(\Gamma, K)\le {{{\rm vol}(\Gamma^*+K)}\over {{\rm vol}(K)}}-1.$$

\medskip\noindent
{\bf Proof.} If $K+X$ is a local packing with $m(\Gamma, K)$ neighbours and satisfying (4.3), then we have
$$K+X\subseteq \Gamma^*+K.$$
Since the interiors of the translates are pairwise disjoint, as it was observed by Hadwiger \cite{had57}, it follows that
$${\rm card}\{ X\}\cdot {\rm vol} (K) \le {\rm vol} (\Gamma^*+K)$$
and therefore
$$m(\Gamma, K)\le {{{\rm vol}(\Gamma^*+K)}\over {{\rm vol}(K)}}-1.$$
The theorem is proved. \hfill{$\Box$}

\medskip\noindent
{\bf Corollary 4.1.} {\it In the octahedral case, we take $r_1=\sqrt{2033}/57$ and $\Gamma_1 =r_1B$. It is easy to see that ${\rm vol} (O)={4\over 3}$ and $\Gamma^*_1=2r_1B$.
By Steiner's formula, one can deduce that
$${\rm vol} (\Gamma_1^*+O)={\rm vol} (O)+4\sqrt{3}\cdot 2r_1+6\sqrt{2}\arccos (1/3)\cdot (2r_1)^2+{\rm vol} (B)\cdot (2r_1)^3.$$
Then, by Theorem $4.2$ one can obtain}
$$m(\Gamma_1, O)\le 40.$$

\medskip\noindent
{\bf Corollary 4.2.} {\it In the cuboctahedral case, we take $r_2=\sqrt{830}/21$ and $\Gamma_2 =r_2B$. It is easy to see that ${\rm vol} (C)={{20}\over 3}$ and $\Gamma^*_2=2r_2B$.
By Steiner's formula, one can deduce that
$${\rm vol} (\Gamma_2^*+C)={\rm vol} (C)+(12+4\sqrt{3})\cdot 2r_2+12\sqrt{2}\arccos \sqrt{1/3}\cdot (2r_2)^2+{\rm vol} (B)\cdot (2r_2)^3.$$
Then, by Theorem $4.2$ it can be shown that}
$$m(\Gamma_2, C)\le 39.$$

\medskip\noindent
{\bf Remark 4.1.} Comparing with Corollary 1.1, the above upper bounds for $m(\Gamma_1, O)$ and $m(\Gamma_2, C)$ are much more hopeful for computer. In fact, by studying $\Gamma_1^*\setminus {\rm int}(2O)$ and
$\Gamma_2^*\setminus {\rm int}(2C)$ in detail these bounds can be further improved.

\vspace{0.6cm}\noindent
{\bf 4.3. The Octahedral Case}

\medskip
\noindent
In this section we take $r_1=\sqrt{2033}/57$, $\Gamma_1=r_1B$ and therefore $\Gamma_1^*=2r_1B$. Assume that $O+X$ is a local packing such that
$$X\setminus \{{\bf o}\}\subset \Gamma_1^*\setminus {\rm int}(2O),$$
where $X=\{ {\bf o}, {\bf x}_1, {\bf x}_2, \ldots , {\bf x}_{m(\Gamma_1, O)}\}$.

Let $r$ be a positive number satisfying ${2\over \sqrt{3}}\le r\le 2r_1=2\sqrt{2033}/57$. It is easy to see that the local packing $O+X$ reduces a spherical packing $(O+{\bf x}_i)\cap \partial (rB)$,
${\bf x}_i\in X$, on $\partial (rB)$. Therefore, detailed study of such two-dimensional packings may provide improvement on estimating $m(\Gamma_1, O)$.

Let $\nu(\cdot )$ denote the area measure on $\partial (rB)$ and define
$$\mu_1(r)=\min_{{\bf x}\in \Gamma_1^*\setminus {\rm int}(2O)}\nu\bigl(\partial (rB)\cap (O+{\bf x})\bigr).$$
Clearly, we have
$$m(\Gamma_1, O)\le {{\nu\bigl(\partial (rB)\bigr)}\over {\mu_1(r)}}={{4\pi r^2}\over {\mu_1(r)}}.\eqno(4.4)$$

\medskip
Next we recall a classic result about spherical cap packing.

\medskip\noindent
{\bf Lemma 4.1 (Moln\'ar \cite{moln52}).} {\it Let $\delta_n$ denote the density of the densest cap packings on $\partial (B)$ with $n$ congruent caps. When $n\ge 3$, we have}
$$\delta_n< {\pi \over {\sqrt{12}}}.$$

\smallskip
Then, Corollary 4.1 can be improved by the following result.

\smallskip\noindent
{\bf Theorem 4.3.}
$$m(\Gamma_1, O)\le 26.$$

\smallskip\noindent
{\bf Proof.} It is well-known that ${{\sqrt{3}}\over 3}B\subset O.$ Based on experiments and observations, we take
$$\rho_1=\sqrt{{1\over 3}+{4\over {57}}\sqrt{{2033}\over 3}}.$$

It is routine to show that
$$\nu\bigl(\partial (\rho_1B)\bigr)=4\pi \left({1\over 3}+{4\over {57}}\sqrt{{2033}\over 3}\right)\eqno(4.5)$$
and
$$\partial (\rho_1B)\cap \left({1\over \sqrt{3}}B+{\bf v}_1\right) =\partial (\rho_1B)\cap \left({1\over \sqrt{3}}B+{\bf v}_2\right),$$
where ${\bf v}_1=\left({2\over \sqrt{3}}, 0, 0\right)$ and ${\bf v}_2=\left(2r_1, 0, 0\right)$. Consequently, for all ${\bf x}\in \Gamma_1^*\setminus {\rm int}(2O),$ we have
$$\nu\left(\partial (\rho_1B)\cap \left({1\over \sqrt{3}}B+{\bf x}\right)\right) \ge \nu\left(\partial (\rho_1B)\cap \left({1\over \sqrt{3}}B+{\bf v}_1\right)\right).$$
Therefore, for any point ${\bf x}\in \Gamma_1^*\setminus {\rm int}(2O)$, the spherical region $\partial (\rho_1B)\cap (O+{\bf x})$ contains a cap which is congruent to $\partial (\rho_1B)\cap \left({1\over \sqrt{3}}B+{\bf v}_1\right)$. Consequently, we get
\begin{align*}
\mu_1(\rho_1)&\ge  \nu\left(\partial (\rho_1B)\cap \left({1\over \sqrt{3}}B+{\bf v}_1\right)\right)\\
&=2\pi \sqrt{{1\over 3}+{4\over {57}}\sqrt{{2033}\over 3}} \left( \sqrt{{1\over 3}+{4\over {57}}\sqrt{{2033}\over 3}}-{{\sqrt{2033}}\over {57}}-{1\over \sqrt{3}}\right).\tag{4.6}
\end{align*}

It is clear that the local packing $O+X$ reduces a spherical cap packing $\left({1\over \sqrt{3}}B+{\bf x}_i\right)\cap \partial (\rho_1B)$,
${\bf x}_i\in X$, on $\partial (\rho_1B)$, which contains a packing with $m(\Gamma_1, O)$ caps all congruent to $\partial (\rho_1B)\cap \left({1\over \sqrt{3}}B+{\bf v}_1\right)$. Thus, by Lemma 4.1, (4.5) and (4.6) we obtain
$$ {{m(\Gamma_1, O) \cdot \mu_1(\rho_1)}\over {\nu \bigl(\partial (\rho_1B)\bigr)}}<{\pi \over \sqrt{12}}$$
and
$$m(\Gamma_1, O)< {\pi \over \sqrt{12}}\cdot {{\nu \bigl(\partial (\rho_1B)\bigr)}\over {\mu_1(\rho_1)}}\le 26.300\ldots.$$
Since $m(\Gamma_1, O)$ is an integer, we have
$$m(\Gamma_1, O)\le 26.$$
The Theorem is proved.\hfill{$\Box$}

\medskip\noindent
{\bf Remark 4.2.} It was shown by Larman and Zong \cite{larm99} that the translative kissing number of an octahedron is $18$. Therefore, the upper bound $26$ is pretty close to the optimal.

\vspace{0.6cm}\noindent
{\bf 4.4. The Cuboctahedral Case}

\medskip
\noindent
In this section we take $r_2=\sqrt{830}/21$, $\Gamma_2=r_2B$ and therefore $\Gamma_2^*=2r_2B$. Assume that $C+X$ is a local packing such that
$$X\setminus \{{\bf o}\}\subset \Gamma_2^*\setminus {\rm int}(2C),$$
where $X=\{ {\bf o}, {\bf x}_1, {\bf x}_2, \ldots , {\bf x}_{m(\Gamma_2, C)}\}$.

Let $r$ be a positive number satisfying $2\le r\le 2r_2=2\sqrt{830}/21$. It is easy to see that the local packing $C+X$ reduces a spherical packing $(C+{\bf x}_i)\cap \partial (rB)$,
${\bf x}_i\in X$, on $\partial (rB)$. Therefore, detailed study of such two-dimensional packings may provide improvement on estimating $m(\Gamma_2, C)$.

Similar to the octahedral case, we define
$$\mu_2(r)=\min_{{\bf x}\in \Gamma_2^*\setminus {\rm int}(2C)}\nu\bigl(\partial (rB)\cap (C+{\bf x})\bigr).$$
Then one can deduce that
$$m(\Gamma_2, C)\le {{\nu\bigl(\partial (rB)\bigr)}\over {\mu_2(r)}}={{4\pi r^2}\over {\mu_2(r)}}.\eqno(4.7)$$

\medskip
Now, Corollary 4.2 can be improved by the following result.

\smallskip\noindent
{\bf Theorem 4.4.}
$$m(\Gamma_2, C)\le 26.$$

\smallskip\noindent
{\bf Proof.} Clearly, we have $B\subset C.$ Based on experiments and observations, we take
$$\rho_2=\sqrt{1+{{4\sqrt{830}}\over {21}}}.$$

It is routine to show that
$$\nu \bigl(\partial (\rho_2B)\bigr)=4\pi \left(1+{{4\sqrt{830}}\over {21}}\right)\eqno(4.8)$$
and
$$\partial (\rho_2B)\cap \left(B+{\bf v}_1\right) =\partial (\rho_2B)\cap \left(B+{\bf v}_2\right),$$
where ${\bf v}_1=(2, 0, 0)$ and ${\bf v}_2=2r_2(1, 0, 0)$. Consequently, for all ${\bf x}\in \Gamma_2^*\setminus {\rm int}(2C),$ we have
$$\nu\bigl(\partial (\rho_2B)\cap \left(B+{\bf x}\right)\bigr) \ge \nu\bigl(\partial (\rho_2B)\cap \left(B+{\bf v}_1\right)\bigr).$$
Therefore, for any point ${\bf x}\in \Gamma_2^*\setminus {\rm int}(2C)$, the spherical region $\partial (\rho_2B)\cap (C+{\bf x})$ contains a cap which is congruent to $\partial (\rho_2B)\cap \left(B+{\bf v}_1\right)$. Consequently, we get
\begin{align*}
\mu_2(\rho_2)&\ge  \nu\bigl(\partial (\rho_2B)\cap \left(B+{\bf v}_1\right)\bigr)\\
&=2\pi \sqrt{1+{{4\sqrt{830}}\over {21}}} \left( \sqrt{1+{{4\sqrt{830}}\over {21}}}-{{\sqrt{830}}\over {21}}-1\right).\tag{4.9}
\end{align*}

It is clear that the local packing $C+X$ reduces a spherical cap packing $(B+{\bf x}_i)\cap \partial (\rho_2B)$,
${\bf x}_i\in X$, on $\partial (\rho_2B)$, which contains a packing with $m(\Gamma_2, C)$ caps all congruent to $\partial (\rho_2B)\cap \left(B+{\bf v}_1\right)$. Thus, by Lemma 4.1, (4.8) and (4.9) we obtain
$$ {{m(\Gamma_2, C) \cdot \mu_2(\rho_2)}\over {\nu \bigl(\partial (\rho_2B)\bigr)}}<{\pi \over \sqrt{12}}$$
and
$$m(\Gamma_2, C)< {\pi \over \sqrt{12}}\cdot {{\nu \bigl(\partial (\rho_2B)\bigr)}\over {\mu_2(\rho_2)}}\le 26.3723\ldots.$$
Since $m(\Gamma_2, C)$ is an integer, we have
$$m(\Gamma_2, C)\le 26.$$
The Theorem is proved.\hfill{$\Box$}

\medskip\noindent
{\bf Remark 4.3.} It was shown by Talata \cite{tala99} that the translative kissing number of a cuboctahedron is $18$. Therefore, the upper bound $26$ is pretty close to the optimal.

\vspace{0.6cm}\noindent
{\bf 4.5. An Application of the Gauss-Bonnet Theorem}

\medskip
\noindent
Let $\Omega $ be a simple region on the surface of $rB$ bounded successively by $n$ circular arcs $C_1$, $C_2$, $\ldots $, $C_n$. Let $\ell_i$ denote the length of $C_i$, let $g_i$ denote the {\it geodesic curvature} of $C_i$, and let $\theta_i$ denote the angle between $C_i$ and $C_{i+1}$ at their meeting point, where $C_{n+1}=C_1$. As a special case of the {\it Gauss-Bonnet Theorem} we have the following formula to compute the area of $\Omega$.

\medskip
\noindent
{\bf Lemma 4.2.}

$$\nu (\Omega)=r^2\left(\sum_{i=1}^n\theta_i+(2-n)\pi -\sum_{i=1}^ng_i\ell_i\right).$$

\medskip
\noindent
{\bf Remark 4.4.} Let $H_i$ denote the hyperplane containing $C_i$ and let $d_i$ denote the distance from the origin ${\bf o}$ to $H_i$. Clearly, $H_i$ divides $\partial (rB)$ into two caps, usually one is bigger than the other.
Then we have
$$g_i=\left\{
\begin{array}{ll}
{{d_i}\over {\sqrt{r^2-d_i^2}}}, & \mbox{if $\Omega $ is in the small cap;}\\
\vspace{-0.35cm}
&\\
{{-d_i}\over {\sqrt{r^2-d_i^2}}}, & \mbox{if $\Omega $ is in the big cap.}
\end{array}\right.
$$
When the two caps are equal, $C_i$ is a part of a great circle on $\partial (rB)$ and therefore $g_i=d_i=0.$

\medskip
To improve the upper bound for $m(\Gamma_1, O)$, we try to determine the value of $\mu_1(\rho_1)$. For this purpose, as shown by Figure 4.3, we take

$$\begin{array}{ll}
{\bf x}_1=\left( {2\over 3}, {2\over 3}, {2\over 3} \right), & {\bf x}'_1=\left( {2\over {57}}\sqrt{{{2033}\over 3}}, {2\over {57}}\sqrt{{{2033}\over 3}}, {2\over {57}}\sqrt{{{2033}\over 3}} \right),\\
\vspace{-0.3cm}
&\\
{\bf x}_2=\left( {2\over 3}-{{10}\over {3\sqrt{57}}}, {2\over 3}-{{10}\over {3\sqrt{57}}}, {2\over 3}+{{20}\over {3\sqrt{57}}}  \right), &
{\bf x}_3=\left( 1-{\sqrt{817}\over {57}}, 0, 1+{\sqrt{817}\over {57}} \right), \\
\vspace{-0.3cm}
&\\
{\bf x}_4= (1, 0, 1),& {\bf x}'_4=\left( {{\sqrt{4066}}\over {57}}, 0, {{\sqrt{4066}}\over {57}} \right).
\end{array}$$

\begin{figure}[ht]
\centering
\includegraphics[height=5cm,width=5cm,angle=0]{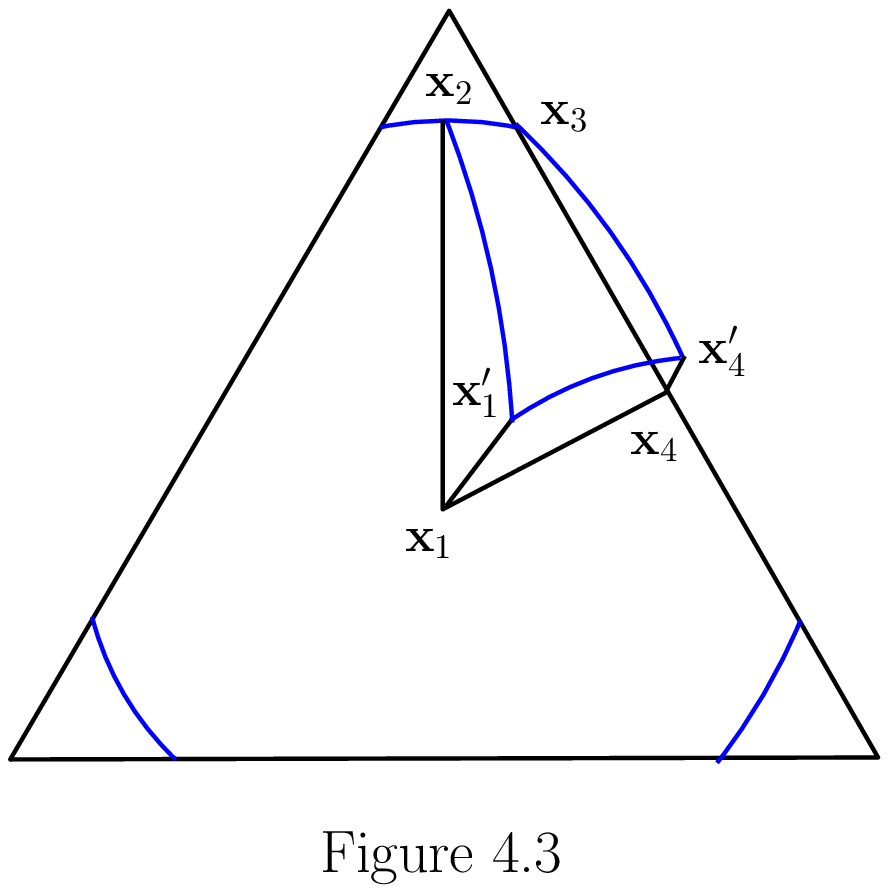}
\end{figure}

\noindent
We define $P_1$ to be the hyperplane passing ${\bf o}$, ${\bf x}_1$ and ${\bf x}_2$, define $P_2$ to be the hyperplane passing ${\bf o}$, ${\bf x}_3$ and ${\bf x}_4$, and define $P_3$ to be the hyperplane passing
${\bf o}$, ${\bf x}_1$ and ${\bf x}_4.$ Let $Q_1$ to be the halfspace bounded by $P_1$ and containing ${\bf x}_4$, let $Q_2$ to be the halfspace bounded by $P_2$ and containing ${\bf x}_1$,
and let $Q_3$ to be the halfspace bounded by $P_3$ and containing ${\bf x}_2$. Finally, we define
$$F_1=\bigl(\Gamma_1^*\setminus {\rm int}(2O)\bigr) \cap Q_1\cap Q_2\cap Q_3.$$

By symmetry, recalling that
$$\rho_1=\sqrt{{1\over 3}+{4\over {57}}\sqrt{{2033}\over 3}},$$
it is easy to show that
$$\mu_1(\rho_1)=\min_{{\bf x}\in F_1} \nu \bigl(\partial (\rho_1B)\cap (O+{\bf x})\bigr).\eqno (4.10)$$
By considering different cases with respect to the shapes of $\partial (\rho_1B)\cap (O+{\bf x})$ and applying Lemma 4.2, one can deduce the following result.

\medskip\noindent
{\bf Lemma 4.3.}
$$\mu_1(\rho_1)= \nu \bigl(\partial (\rho_1B)\cap (O+{\bf x}_2)\bigr).$$

\medskip
As a direct consequence of (4.4) and Lemma 4.3 we obtain the following improvement of Theorem 4.3.

\medskip\noindent
{\bf Theorem 4.5.}
$$m(\Gamma_1, O)\le 22.$$

\bigskip
To improve the upper bound for $m(\Gamma_2, C)$, we try to determine the value of $\mu_2(\rho_2)$. For this purpose, as shown by Figure 4.4, we take

$$\begin{array}{ll}
{\bf y}_1=\left( 2, 0, 0 \right), & {\bf y}'_1=\left( {{2\sqrt{830}}\over {21}}, 0, 0  \right), \\
\vspace{-0.3cm}
&\\
{\bf y}_2=\left( 2, 1, 1 \right), & {\bf y}'_2=\left( {4\over {21}}\sqrt{{415}\over 3}, {2\over {21}}\sqrt{{415}\over 3}, {2\over {21}}\sqrt{{415}\over 3} \right), \\
\vspace{-0.3cm}
&\\
{\bf y}_3=\left( {4\over 3}, {4\over 3}, {4\over 3} \right), & {\bf y}'_3=\left( {2\over {21}}\sqrt{{830}\over 3}, {2\over {21}}\sqrt{{830}\over 3}, {2\over {21}}\sqrt{{830}\over 3} \right), \\
\vspace{-0.3cm}
&\\
{\bf y}_4=\left( {4\over 3}+{{22\sqrt{3}}\over {63}}, {4\over 3}-{{44\sqrt{3}}\over {63}}, {4\over 3}+{{22\sqrt{3}}\over {63}} \right),& {\bf y}_5=\left( 2, 1-{\sqrt{337}\over {21}}, 1+{\sqrt{337}\over {21}} \right),\\
\vspace{-0.3cm}
&\\
{\bf y}_6=\left(2, 0, {{2\sqrt{389}}\over {21}}\right). &
\end{array}$$

\begin{figure}[ht]
\centering
\includegraphics[height=5.5cm,width=5.5cm,angle=0]{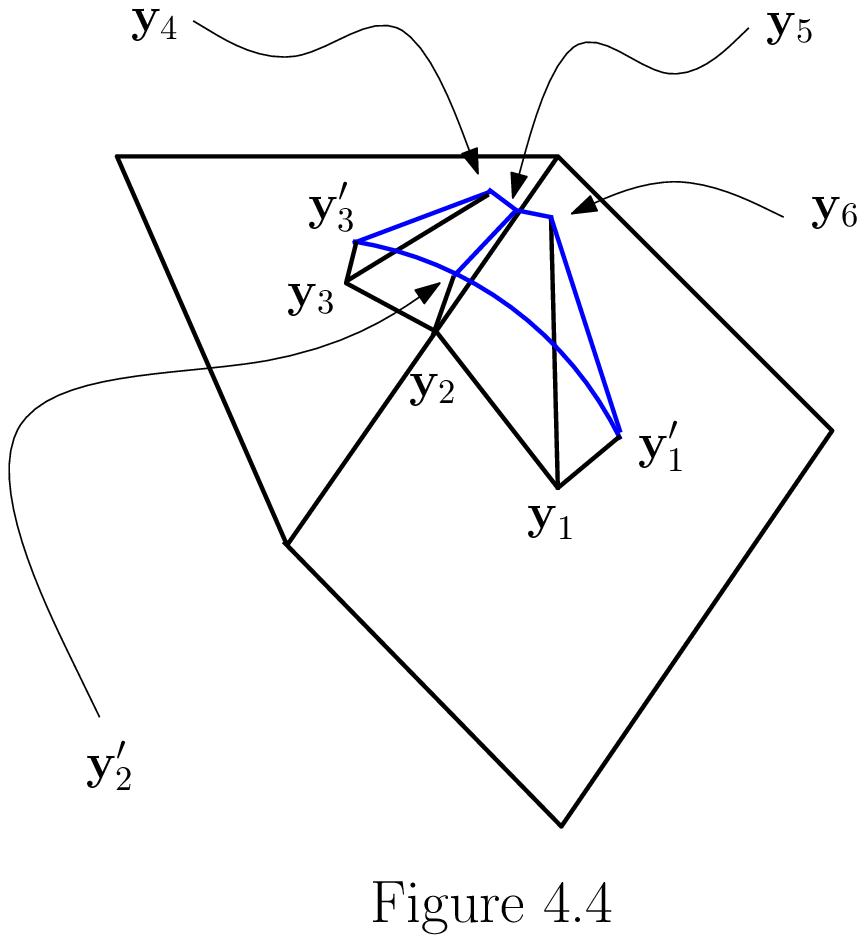}
\end{figure}

\noindent
We define $P_1$ to be the hyperplane passing ${\bf o}$, ${\bf y}_1$ and ${\bf y}_2$, define $P_2$ to be the hyperplane passing ${\bf o}$, ${\bf y}_3$ and ${\bf y}_4$, and define $P_3$ to be the hyperplane passing
${\bf o}$, ${\bf y}_1$ and ${\bf y}_4.$ Let $Q_1$ to be the halfspace bounded by $P_1$ and containing ${\bf y}_4$, let $Q_2$ to be the halfspace bounded by $P_2$ and containing ${\bf y}_1$,
and let $Q_3$ to be the halfspace bounded by $P_3$ and containing ${\bf y}_2$. Finally, we define
$$F_2=\bigl(\Gamma_2^*\setminus {\rm int}(2C)\bigr) \cap Q_1\cap Q_2\cap Q_3.$$

By symmetry, recalling that
$$\rho_2=\sqrt{1+{{4\sqrt{830}}\over {21}}},$$
it is easy to show that
$$\mu_2(\rho_2)=\min_{{\bf x}\in F_2} \nu \bigl(\partial (\rho_2B)\cap (C+{\bf x})\bigr).\eqno (4.11)$$
Similar to the octahedral case, by considering different cases with respect to the shapes of $\partial (\rho_2B)\cap (C+{\bf x})$ and applying Lemma 4.2, one can deduce the following result.

\medskip\noindent
{\bf Lemma 4.4.}
$$\mu_2(\rho_2)= \nu \bigl(\partial (\rho_2B)\cap (C+{\bf y}'_{3})\bigr).$$

\medskip\noindent
{\bf Theorem 4.6.}
$$m(\Gamma_2, C)\le 22.$$

\vspace{1cm}
\centerline{\Large\bf 5. Conclusion}

\vspace{0.6cm}\noindent
If Conjecture Z is true, it can be proved by dealing with optimization problems corresponding to triangulated color graphs with at most 22 vertices and four or seven colors, respectively. Of course, the target functions are more complicated than (2.6) for the reason of the truncater.

\vspace{0.8cm}\noindent
{\large\bf Acknowledgement.} This work is supported by 973 Program 2013CB834201 and the Chang Jiang Scholars Program of China. For helpful discussions, I am grateful to professor R. J. Gardner and professor F. Santos.

\bigskip
\bibliographystyle{amsplain}

\bigskip\noindent
Center for Applied Mathematics, Tianjin University, Tianjin 300072, China

\smallskip\noindent
cmzong@math.pku.edu.cn

\end{document}